\documentclass[reqno]{amsart}

\usepackage{amsthm,amsmath,amsfonts,amssymb,epsfig,graphicx,latexsym,float}

\newtheorem{definition}{Definition}

\newtheorem{remark}[definition]{Remark}
\newtheorem{assumption}[definition]{Assumption}

\newtheorem{example}[definition]{Example}
\newtheorem{notation}[definition]{Notation}

\renewcommand{\theequation}{\arabic{section}.\arabic{equation}}
\renewcommand{\thetable}{\arabic{section}.\arabic{table}}
\renewcommand{\thefigure}{\arabic{section}.\arabic{figure}}

\begin{document}

\title[Low-Mach-number--slenderness limit for elastic Cosserat rods]
{Low-Mach-number--slenderness limit \\ for elastic Cosserat rods}

\author[F.~Baus]{Franziska Baus$^{1}$}
\author[A.~Klar]{Axel Klar$^{1}$}
\author[N.~Marheineke]{Nicole Marheineke$^{2,\star}$}
\author[R.~Wegener]{Raimund Wegener$^{3}$} 

\date{\today\\
$^\star$ \textit{Corresponding author}, email: marheineke@math.fau.de, phone: +49\,9131\,85\,67214, fax: +49\,9131\,85\,67225\\
$^1$ TU Kaiserslautern, Fachbereich Mathematik, Erwin-Schr\"odinger-Str.~48, D-67663 Kaiserslautern, Germany\\
$^2$ FAU Erlangen-N\"urnberg, Lehrstuhl Angewandte Mathematik I, Cauerstr.~11, D-91058 Erlangen, Germany\\
$^3$ Fraunhofer ITWM, Fraunhofer Platz 1, D-67663 Kaiserslautern, Germany
}

\begin{abstract}
This paper deals with the relation of the dynamic elastic Cosserat rod model and the Kirchhoff beam equations. We show that the Kirchhoff beam without angular inertia is the asymptotic limit of the Cosserat rod, as the slenderness parameter (ratio between rod diameter and length) and the Mach number (ratio between rod velocity and typical speed of sound) approach zero, i.e.\ low-Mach-number--slenderness limit. 
The asymptotic framework is exact up to fourth order in the small parameter and reveals a mathematical structure that allows a uniform handling of the transition regime between the models. To investigate this regime numerically, we apply a scheme that is based on a spatial Gauss-Legendre collocation and an $\alpha$-method in time.
\end{abstract}

\maketitle

{\sc AMS-Classification.} 74K10, 35Q74, 35B40, 35C20, 65Mxx
%\keywords{{\sc Keywords.} dynamic elastic Cosserat rod, Kirchhoff beam, low-Mach-number--slenderness limit, asymptotic analysis, asymptotic-preserving scheme}
%\subjclass{{\sc AMS-Classification.} 74K10, 35Q74, 35B49, 35C20, 65Mxx}

%%%%%%%%%%%%%%%%%%%%%%%%%%%%%%%%%%%%%%%%%%%%%%%%%%%%%%%%%%%%%%%%%%
\section{Introduction}\label{sec:1}
The elastic rod theory is an old, extensively studied, but still current topic of research. Its foundations date back to, among others, Bernoulli, Kirchhoff \cite{kirchhoff:p:1859}, the brothers Cosserat \cite{cosserat:b:1909} and Love \cite{love:b:1927}. A comprehensive overview from today's perspective is given in, for example \cite{antman:b:2006, rubin:b:2000}. The investigations range from analytical aspects such as solution theory, stability and Hamiltonian structure to numerical methods (geometrically exact approaches, exact energy-momentum conserving algorithms and symplectic schemes), of the broad literature see e.g.\  \cite{maddocks:p:1984, maddocks:p:1994, simo:p:1988, chouaieb:p:2004, coleman:p:1993, langer:p:1996} and \cite{simo:p:1985, simo:p:1986, simo:p:1995, romero:p:2002, bergou:p:2008, audoly:b:2010}. Equally large is the field of applications: engineering mechanics (truss, fiber-reinforced materials, non-woven textiles, paper), biomolecular science (DNA, bacterial fibers), computer graphics etc.

In this paper the specific focus lies on the asymptotic investigation of the relation between the dynamic Cosserat rod model and the Kirchhoff beam equations to describe the motion of slender elastic bodies. A rod in the special Cosserat theory \cite{antman:b:2006} is represented by two constitutive elements: a parameterized time-dependent curve $\mathbf{r}$ and an attached orthonormal director triad $\{\mathbf{d_1},\mathbf{d_2},\mathbf{d_3}\}$ that specify the position and the orientation of the material cross-sections. Its deformation is due to tension, shear, bending and torsion. For general elastic material laws we derive beam equations of Kirchhoff-type in an asymptotic analysis, as the slenderness parameter $\epsilon$ (ratio between rod diameter and length) and the Mach number $\mathrm{Ma}$ (ratio between rod velocity and typical speed of sound) approach zero. In the combined asymptotic limit $\epsilon \rightarrow 0$, $\mathrm{Ma}\rightarrow 0$ with $\mathrm{Ma}/\epsilon =\mu = \mathrm{const}$ to which we refer as low-Mach-number--slenderness limit, the model equations show two characteristic changes: 1) the contact force becomes a variable to a constraint on the kinematics and the respective material law decouples from the  model;  2) the angular inertia terms vanish.
The terminology Kirchhoff stands for an inextensible and unshearable rod. Classically, the Kirchhoff constraint relates the beam tangent and the director triad, it is $\partial_s\mathbf{r}=\mathbf{d_3}$, see e.g.\ \cite{antman:b:2006, coleman:p:1993, langer:p:1996}. As consequence of the constraint, the contact force and the strain change their roles in the model system. Their proportionality operator is determined by the material law for the contact force, it can be identified with the moduli of linear elasticity theory. But in contrast to the linear theory whose validity is restricted to small deformations, large deformations are allowed as the material law can be chosen arbitrarily and the proportionality operator is just its linearization.  
Since namings of different beam models is frequently problematic and inconsistent in literature, we point out that we call the limit model a Kirchhoff beam when the system is equipped with the Euler-Bernoulli relation for the contact couple. In that case the deduced limit model that has additionally no angular inertia terms is also known as Kirchhoff-Love equations \cite{landau:b:1970}. 

Deriving beam models with different degrees of freedom from the three-dimensional theory of elasticity is topic in several works, see for example the asymptotic limits in \cite{coleman:p:1993, coleman:p:1995}. The low-Mach-number--slenderness asymptotics presented in this paper reveals a direct scaling between the Cosserat rod and the Kirchhoff beam without angular inertia and complements the previous works. The result has its analogue in the fluid dynamics where the low-Mach-number asymptotics describes the transition from a compressible to an incompressible fluid \cite{lions:p:1998, alazard:p:2006, meister:p:1999}. For elastic rods the speed of sound is not a unique quantity, instead compression and shear waves induce different speeds of sound that imply different Mach numbers. For the asymptotic derivation, it is assumed that there exists a typical magnitude of the contact force. Such a scaling presupposes that the Mach numbers associated with compression and shear behave similar. This corresponds particularly to a constant Poisson ratio for a linear isotropic material, as it is often considered in macroscopic applications.

The asymptotic rod behavior is here characterized by the slenderness parameter $\epsilon$ and the (typical) Mach number $\mathrm{Ma}$ which we relate according to $\mathrm{Ma}/\epsilon=\mu=\mathrm{const}$. In the transition regime for small $\epsilon$ ($\epsilon \rightarrow 0$) the type of the model equations changes: time derivatives degenerate, the system becomes stiff with a constraint. The asymptotic framework reveals a mathematical structure that allows a uniform robust numerical handling of this regime. We show that the asymptotic systems (limit system and its first-order correction) which follow from an even power series expansion are together exact up to order $\mathcal{O}(\epsilon^4)$. Moreover, the conservation of energy is ensured, if hyper-elastic material laws, external potential forces and appropriate boundary conditions are presupposed. Having a Newton method for the numerical treatment in mind, solving both asymptotic systems is of similar computational effort as solving the original $\epsilon$-dependent system while it is much more accurate and robust in computing the influence of small $\epsilon$-values. We point out that the underlying discretization is replaceable. In this paper we propose an asymptotic-preserving scheme in the spirit of the generalized $\alpha$-method that can be switched between energy-conserving and dissipative \cite{chung:p:1993, sobottka:p:2008}. This is advantageous for applications with external non-potential forces, such as fiber-fluid interactions in non-woven manufacturing \cite{marheineke:p:2011, klar:p:2009} or hair simulations \cite{bergou:p:2008,sobottka:p:2008}. We refer to \cite{simo:p:1995, romero:p:2002} for exact energy-momentum conserving algorithms and to \cite{bergou:p:2008, audoly:b:2010} for discrete geometric approaches. Moreover, we refer to the research work in fluid dynamics where the low-Mach-number asymptotics has been used to develop and extend numerical schemes for the compressible-incompressible transition regime, e.g.\ \cite{klein:p:1995,guillard:p:1999}.

This paper is structured as follows. Starting with a brief introduction into the special Cosserat rod theory, we present the low-Mach-number--slenderness asymptotics in Section~\ref{sec:2}. We derive the beam equations of Kirchhoff-type for general elastic material laws. To study numerically the performance of the asymptotic framework in the transition regime for small $\epsilon$, we consider a two-dimensional Euler-Bernoulli cantilever beam in Section~\ref{sec:3}. Three variants of the underlying rod model can be distinguished, depending on the chosen kinematic/geometric formulation. They imply temporal or spatial differential-algebraic systems of different index after semi-discretization in space or time, respectively. In the Appendix we provide details to the underlying numerical scheme that is applicable to all formulations in the transition regime. It is based on a Gauss-Legendre collocation in space (finite differences) and an $\alpha$-method in time, but can be also viewed as a conservative finite-volume method.

%%%%%%%%%%%%%%%%%%%%%%%%%%%%%%%%%%%%%%%%%%%%%
\setcounter{equation}{0} \setcounter{figure}{0}
\section{Asymptotic relation between elastic Cosserat rod and Kirchhoff beam}\label{sec:2}

\subsection{Special Cosserat rod theory}

An elastic thread is a slender long body, i.e.\ a rod in three-dimensional continuum mechanics. Because of its geometry the dynamics might be reduced to a one-dimensional description by averaging the underlying balance laws over its cross-sections. This procedure is based on the assumption that the displacement field in each cross-section can be expressed in terms of a finite number of vector- and tensor-valued quantities. The special Cosserat rod theory \cite{antman:b:2006} consists of only two constitutive elements in the three-dimensional Euclidean space $\mathbb{E}^3$, a curve $\mathbf{r}:\mathcal{D}\rightarrow\mathbb{E}^3$ specifying the position and an orthonormal director triad $\{\mathbf{d_1},\mathbf{d_2},\mathbf{d_3}\}:\mathcal{D}\rightarrow\mathbb{E}^3$ characterizing the orientation of the cross-sections. In $\mathcal{D}=\{(s,t)\in \mathbb{R}^2 \,|\, s \in [s_a,s_b],\, t> 0\}$, the parameter $s$ denotes the material cross-section (material point) and $t$ the time. The rod system involves four kinematic/geometric and two dynamic equations (balance laws for linear and angular momentum)
\begin{equation}\label{eq:rod}
\begin{aligned}
\partial_t \mathbf{r} = \mathbf {v},  \hspace*{1.cm} & \hspace*{1.5cm} \partial_s \mathbf{r} = \boldsymbol{\tau}\\
 \partial_t \mathbf{d_k} = \boldsymbol{\omega} \times \mathbf{d_k}, \hspace*{1.cm}&
\hspace*{1.5cm} \partial_s \mathbf{d_k}=\boldsymbol{\kappa}\times \mathbf{d_k}, \qquad k=1,2,3\\ 
\partial_t((\rho A)  \mathbf{v}) &= \partial_s \mathbf{n} + \mathbf{f}\\
\partial_t ((\rho \mathbf{J})\cdot \boldsymbol{\omega}) &= \partial_s \mathbf{m} + \boldsymbol{\tau}\times\mathbf{n} + \mathbf{l}
\end{aligned}
\end{equation}
with tangent $\boldsymbol{\tau}$, generalized curvature $\boldsymbol{\kappa}$, linear velocity $\mathbf{v}$ and angular velocity $\boldsymbol{\omega}$.  The mass line density $(\rho A)$ is time-independent for materially closed systems, but modeled as time-dependent for applications, such as evaporation and aggregation. The tensor-valued moment of inertia of the cross-sections $(\rho \mathbf{J})$ depends on the configuration of the triad and is hence always time-dependent. To close the system of equations we need to specify the external loads $\mathbf{f}$ and $\mathbf{l}$, boundary and initial conditions as well as elastic material laws for the contact force $\mathbf{n}$ and couple $\mathbf{m}$. In this paper we consider time-independent mass properties of the rod, i.e.\ $\partial_t(\rho A)=0$ and $\partial_t(\mathbf{d_i}\cdot (\rho \mathbf{J})\cdot \mathbf{d_j})=0$ for $i,j=1,2,3$,  and neglect external couples, i.e.\ $\mathbf{l}=\mathbf{0}$, for reasons of a simple presentation. However, extensions are straightforward possible. 
Alternatively to \eqref{eq:rod}, the rod system can be also set up by using the compatibility relations between the kinematics and geometry for the curve and for the triad. These two compatibility conditions 
\begin{align*}
\partial_t \boldsymbol{\tau} &= \partial_s \mathbf{v}\\
\partial_t \boldsymbol{\kappa} &= \partial_s \boldsymbol{\omega} + \boldsymbol{\omega}\times \boldsymbol{\kappa}
\end{align*}
replace then either the two kinematic equations or the two geometric equations.

\begin{notation}[Model variants (M), (T), (S)]\label{model:rod}Depending on the chosen kinematic/geometric formulation, we distinguish three variants of the rod model: 
\begin{itemize}
\item[(M)] kinematic and geometric equations and balance laws (cf.\ \eqref{eq:rod})
\item[(T)] kinematic equations, compatibility conditions and balance laws
\item[(S)] geometric equations, compatibility conditions and balance laws
\end{itemize}
Throughout the paper we often summarize the three systems in a compact form (M-T-S) for readability. Apart from the balance laws, (M-T-S) contains then all six kinematic, geometric and compatibility equations together, see for the first time in \eqref{eq:R_Cosserat}.
\end{notation}
The model variant (T) yields a hyperbolic system for a dynamic elastic rod, whereas (S) is very suitable for the transition to a stationary consideration. 
Presupposing hyper-elastic constitutive relations and external potential loads, the variant (M) is known as bi-Hamiltonian form in literature \cite{simo:p:1988a, dichmann:p:1996}.  The choice of the formulation can affect the numerical simulations as we will comment on in Section~\ref{sec:3} and the Appendix. However, the distinction plays no role for the asymptotics, thus we make use of the compact notation (M-T-S) in the following.

For the objective formulation of the material laws it is useful to rewrite the rod model in the director basis. To an arbitrary vector field $\mathbf{z}=\sum_{i=1}^3 z_i \mathbf{d_i}= \sum_{i=1}^3 \breve z_i \mathbf{e_i}\in \mathbb{E}^3$, we indicate the coordinate tuples corresponding to the director basis $\{\mathbf{d_1},\mathbf{d_2}, \mathbf{d_3}\}$ and to a fixed outer Cartesian basis $\{\mathbf{e_1},\mathbf{e_2},\mathbf{e_3}\}$ by $\mathsf{z}=(z_1, z_2, z_3)\in \mathbb{R}^3$ and $\mathsf{\breve z}=(\breve z_1, \breve z_2, \breve z_3)\in \mathbb{R}^3$, respectively. The director basis can be transformed into the outer basis by the tensor-valued rotation $\mathbf{D}$, i.e.\ $\mathbf{D}=\mathbf{e_i}\otimes \mathbf{d_i} =D_{ij} \mathbf{e_i}\otimes \mathbf{e_j}\in \mathbb{E}^3\otimes \mathbb{E}^3$ with the associated orthogonal matrix $\mathsf{D}=(D_{ij})=(\mathbf{d_i}\cdot \mathbf{e_j})\in SO(3)$. For the coordinates, the relation $\mathsf{D} \cdot \mathsf{\breve z}=\mathsf{z}$ holds -- as well as $\mathsf{D}\cdot \partial_t \mathsf{\breve z}=\partial_t \mathsf{z}+\mathsf{\omega}\times \mathsf{z}$ and $\mathsf{D}\cdot \partial_s\mathsf{\breve z}=\partial_s \mathsf{z}+\mathsf{\kappa}\times \mathsf{z}$. The rotation matrix $\mathsf{D}$ can be parametrized, for example, in Euler angles or unit quaternions. Moreover, canonical basis vectors in $\mathbb{R}^3$ are denoted by $\mathsf{e}_i$, $i=1,2,3$, e.g.\ $\mathsf{e_1}=(1,0,0)$. In the director basis, the Cosserat rod model has the following form (M-T-S)
\begin{equation}\label{eq:R_Cosserat}
\begin{aligned} 
\mathsf{D} \cdot \partial_t \mathsf{\breve r} = \mathsf{v}, \hspace{1cm}& \hspace{4cm} \mathsf{D} \cdot \partial_s \mathsf{\breve r}=\mathsf{\tau}\\ 
\partial_t \mathsf{D} = -\mathsf{\omega}\times \mathsf{D}, \hspace{1cm}& \hspace{4cm} \partial_s \mathsf{D} = -\mathsf{\kappa}\times \mathsf{D}\\
\partial_t \mathsf{\tau} &= \partial_s \mathsf{v}+\mathsf{\kappa}\times \mathsf{v} - \mathsf{\omega} \times \mathsf{\tau}\\
\partial_t \mathsf{\kappa} &= \partial_s \mathsf{\omega} +\mathsf{\kappa}\times \mathsf{\omega} \\[1.5ex]
(\rho A)\partial_t\mathsf{v}&= \partial_s \mathsf{n} + \mathsf{\kappa}\times \mathsf{n} - \mathsf{\omega}\times (\rho A) \mathsf{ v} +\mathsf{D}\cdot \mathsf{\breve f}\\
(\rho \mathsf{J})\cdot \partial_t \mathsf{\omega}  &= 
\partial_s \mathsf{m} + \mathsf{\kappa}\times \mathsf{m} + \mathsf{\tau}\times
\mathsf{n}  - \mathsf{\omega} \times ((\rho \mathsf{J})\cdot \mathsf{\omega})
\end{aligned}
\end{equation}
with general elastic material laws
\begin{align*}
\mathsf{n}(s,t)=\mathcal{N}(\mathsf{\tau}(s,t),\mathsf{\kappa}(s,t),s), \qquad \qquad \mathsf{m}(s,t)=\mathcal{M}(\mathsf{\tau}(s,t),\mathsf{\kappa}(s,t),s).
\end{align*}
Following the assumption on the mass properties, $(\rho \mathsf{J})$ is here the time-independent representation of the inertia tensor in the director basis. In the stated general form the objective elastic constitutive functions $\mathcal{N}$ and $\mathcal{M}$ depend on the strain variables (coordinate tuples) $\mathsf{\tau}$ and $\mathsf{\kappa}$ whose components measure shear $\tau_1$, $\tau_2$, dilatation $\tau_3$, flexure $\kappa_1$, $\kappa_2$ and torsion $\kappa_3$. 
The subsequent asymptotic analysis is valid for material laws that possess the following properties.

\begin{assumption}[Properties of the material law $\mathcal{N}$]\label{ass:1} Assume that the material law $\mathcal{N}$ for the contact forces fulfills the following conditions:
\begin{itemize}
 \item[a)] $\mathcal{N}(\mathsf{\tau},\mathsf{\kappa},s)=0$ holds if and only if $\mathsf{\tau}=\overset{\scriptscriptstyle \circ}{\mathsf{\tau}}(s)$.
 \item[b)] $\mathcal{N}$ is sufficiently regular and $\partial_\tau \mathcal{N}(\overset{\scriptscriptstyle \circ}{\mathsf{\tau}}(s),\mathsf{\kappa},s)$ is an invertible linear operator (matrix) for arbitrary $\mathsf{\kappa}$ and $s$.
\end{itemize}
\end{assumption}

\begin{remark}
Assumption~\ref{ass:1}a) presumes that the family of stress-free configurations is uniquely determined by a strain field $\overset{\scriptscriptstyle \circ}{\mathsf{\tau}}$. Without loss of generality an arc-length parameterization can be chosen for these configurations, implying $\|\overset{\scriptscriptstyle \circ}{\mathsf{\tau}}\|=1$. 
Assumption~\ref{ass:1}b) is naturally satisfied by the monotonicity condition on the constitutive laws \cite{antman:b:2006}. In case of hyper-elastic constitutive relations, i.e.\ $\mathcal{N}=\partial_\mathsf{\tau}\Psi$ and $\mathcal{M}=\partial_\mathsf{\kappa}\Psi$,
it is expressed in terms of a strictly convex elastic potential $\Psi$. 
\end{remark}

\begin{example}[Specification of material laws]\label{ex:1}
As an example for a (hyper-)elastic constitutive law satisfying the assumptions above one can think of a Timoshenko beam under small deformations that is equipped with the affine linear Euler-Bernoulli relation for the contact couple, i.e.\
\begin{align*}
\mathcal{N}(\mathsf{\tau},\mathsf{\kappa},s)=(\mathsf{EA})(s)\cdot (\mathsf{\tau}-\overset{\scriptscriptstyle \circ}{\mathsf{\tau}}(s)), \qquad \qquad 
\mathcal{M}(\mathsf{\tau},\mathsf{\kappa},s)=(\mathsf{EJ})(s)\cdot (\mathsf{\kappa}-\overset{\scriptscriptstyle \circ}{\mathsf{\kappa}}(s)).
\end{align*}
In literature the Timoshenko beam is in general associated with $\overset{\scriptscriptstyle \circ}{\mathsf{\tau}}(s)=\mathsf{e_3}$.
The positive-definite tensor-valued functions $(\mathsf{EA})$ and $(\mathsf{EJ})$ can be expressed by the scalar-valued properties of the cross-section associated Young's modulus $(EA)$ and the shearing modulus $(GA)$ -- in combination with the mass line density $(\rho A)$ and the moment of inertia $(\rho I)$. In case of a geometry with circular cross-sections, we obtain the diagonal forms\footnote{The diagonal forms of $(\rho \mathsf{J})$ and of the linear operator associated with $\mathcal{M}$ are necessary for the embedding into a 2d test scenario. But note that for this purpose the special choice of circular cross-sections is not compulsory.}
\begin{align*}
(\rho \mathsf{J})=(\rho I)\mathsf{P}_2, \qquad \qquad
(\mathsf{EA})=\frac{(EA)}{a} \mathsf{P}_a, \qquad \qquad (\mathsf{EJ})=\frac{(\rho I) (EA)}{(\rho A)} \mathsf{P}_{2/a}
\end{align*}
where $\mathsf{P}_k=\mathrm{diag}(1,1,k)$, $k\in \mathbb{R}$, and $a=(EA)/(GA)=2(1+\nu)$ with Poisson ratio $\nu\in [0,0.5)$.
The physical quantities $(\rho A)$, $(\rho I)$ and $(EA)$ are particularly constant for an homogeneous thread. 
\end{example}

%%%%%%%%%%%
\subsection{Asymptotic low-Mach-number--slenderness limit}
Proceeding from the Cosserat rod model, we derive beam equations of Kirchhoff-type (a generalized string model) as the combined slenderness and low-Mach-number limits in this subsection.  We consider an elastic thread. It has various physical properties for which we choose typical constant reference values that we mark by the subscript $_\star$, i.e.\ mass line density $(\rho A)_\star$, moment of inertia $(\rho I)_\star$, length $L_\star$, magnitude of mean velocity $V_\star$ and contact force $N_\star$. When making \eqref{eq:R_Cosserat} dimensionless, we can characterize the dynamics of the thread by help of two dimensionless numbers: the slenderness parameter $\epsilon$ that is the ratio between the thread's diameter and its length as well as the (typical) Mach number $\mathrm{Ma}$ that is the ratio between the thread's velocity and the speed of sound
\begin{align*}
 \epsilon=\frac{1}{L_\star}\sqrt{\frac{(\rho I)_\star}{(\rho A)_\star}}, \qquad \qquad \mathrm{Ma}=V_\star\sqrt{\frac{(\rho A)_\star}{N_\star}}.
\end{align*}
The speed of sound is here determined by the typical contact force $N_\star$. In the classical sense of compression and shear waves, an elastic thread has certainly different speeds of sounds. Hence, the applied scaling presupposes that the associated Mach numbers behave similar. This means in the special case of a linear isotropic material (cf.\ Example~\ref{ex:1}) that the typical Young's modulus could be chosen, $N_\star=(EA)_\star$, and that the Poisson ratio satisfies $\nu=\mathrm{const}$ in the asymptotics.

For the scaling, we use the following reference values
\begin{align*}
s_0=r_0=L_\star, \quad v_0=V_\star, \quad t_0=\frac{L_\star}{V_\star}, \quad \omega_0=\frac{V_\star}{L_\star}, \quad \kappa_0=\frac{1}{L_\star}, \quad (\rho A)_0=(\rho A)_\star. \quad (\rho J)_0=(\rho I)_\star\\
n_0=(\rho A)_\star V_\star^2, \quad m_0=(\rho A)_\star L_\star V_\star^2, \quad \mathcal{N}_0=N_\star, \quad \mathcal{M}_0=\frac{(\rho I)_\star N_\star}{(\rho A)_\star L_\star}, \quad f_0=\frac{(\rho A)_\star V_\star^2}{L_\star}\quad
\end{align*}
and introduce to every dimensional variable $\mathsf{z}$ the associated dimensionless one $\bar{\mathsf z}$ as
\begin{align*}
\mathsf{z}(s_0\bar{s},t_0\bar{t})=z_0 \bar{\mathsf{z}}(\bar{s},\bar{t}), \qquad \qquad s={s_0} \bar s, \quad t= {t_0} \bar t.
\end{align*}
Skipping the superscript $\bar{\,\,}$ for readability, the dimensionless rod system is then given by
\begin{equation}\label{eq:dimless_Rod}
\begin{aligned} 
\mathsf{D} \cdot \partial_t \mathsf{\breve r} = \mathsf{v}, \hspace{1cm} &  \hspace{4cm} \mathsf{D} \cdot \partial_s \mathsf{\breve r}=\mathsf{\tau}\\ 
\partial_t \mathsf{D} = -\mathsf{\omega}\times \mathsf{D},  \hspace{1cm} &  \hspace{4cm}  \partial_s \mathsf{D} = -\mathsf{\kappa}\times \mathsf{D}\\
\partial_t \mathsf{\tau} &= \partial_s \mathsf{v}+\mathsf{\kappa}\times \mathsf{v} + \mathsf{\tau}\times \omega\\
\partial_t \mathsf{\kappa} &= \partial_s \mathsf{\omega} +\mathsf{\kappa}\times \mathsf{\omega} \\[1.5ex] 
(\rho A)\partial_t\mathsf{v} &= \partial_s \mathsf{n} + \mathsf{\kappa}\times \mathsf{n} +(\rho A) \mathsf{ v} \times \mathsf{\omega}+\mathsf{D}\cdot \mathsf{\breve f}\\ 
\epsilon^2(\rho \mathsf{J}) \cdot \partial_t \mathsf{\omega}  &= 
\partial_s \mathsf{m} + \mathsf{\kappa}\times \mathsf{m} + \mathsf{\tau}\times
\mathsf{n}  +\epsilon^2((\rho \mathsf{J})\cdot \mathsf{\omega})\times \mathsf{\omega} 
\end{aligned}
\end{equation}
with general elastic material laws ($\mathcal{N}$ satisfying Assumption~\ref{ass:1})
\begin{align*}
\mathsf{n}(s,t)=\frac{1}{\mathrm{Ma}^2}&\mathcal{N}(\mathsf{\tau}(s,t),\mathsf{\kappa}(s,t),s), \qquad \qquad \mathsf{m}(s,t)=\frac{\epsilon^2}{\mathrm{Ma}^2}\mathcal{M}(\mathsf{\tau}(s,t),\mathsf{\kappa}(s,t),s).
\end{align*}

To establish an asymptotic relation between the Cosserat rod \eqref{eq:dimless_Rod} and a beam of Kirchhoff-type, we consider the limits 
\begin{align*}
\epsilon\rightarrow 0, \quad \mathrm{Ma}\rightarrow 0 \quad \quad \text{with} \quad \frac{\mathrm{Ma}}{\epsilon}=\mu, \quad \mu=\text{const}>0
\end{align*}
(i.e.\ asymptotic limit along straight lines in the $(\epsilon, \mathrm{Ma})$-plane with slope $\mu$).
As the small parameter $\epsilon$ appears only in second power in the equations ($\epsilon^2$, $\mathrm{Ma}^2=\mu^2 \epsilon^2$), we expand all quantities in an even power series with respect to the slenderness parameter,  i.e.\ $\mathsf{z}=\sum_{i=0}^\infty \epsilon^{2i}\mathsf{z}^{(i)}$.
The contact force $\mathsf{n}(s,t)=(\mu\epsilon)^{-2}\mathcal{N}(\mathsf{\tau}(s,t),\mathsf{\kappa}(s,t),s)$ implies
\begin{align*}
 0=\mathcal{N}(\mathsf{\tau}^{(0)}(s,t),\mathsf{\kappa}^{(0)}(s,t),s), \quad \text{ as } \epsilon = 0, \quad \quad \text{ hence } \mathsf{\tau}^{(0)}(s,t) = \overset{\scriptscriptstyle \circ}{\mathsf{\tau}}(s)
\end{align*}
according to Assumption~\ref{ass:1}a).
Consequently, we set the strain to be $\mathsf{\tau}=\overset{\scriptscriptstyle \circ}{\mathsf{\tau}}+\epsilon^2 \mathsf{\xi}$ with a function $\mathsf{\xi}\sim \mathcal{O}(1)$. Then the contact force becomes
$\mathsf{n}(s,t)=(\mu\epsilon)^{-2}\mathcal{N}(\overset{\scriptscriptstyle \circ}{\mathsf{\tau}}(s)+\epsilon^2\mathsf{\xi}(s,t),\mathsf{\kappa}(s,t),s)$, whose leading order we obtain via Taylor expansion
\begin{align*}
\mathsf{n}^{(0)}(s,t)=\mu^{-2}\,(\partial_\tau \mathcal{N}(\overset{\scriptscriptstyle \circ}{\mathsf{\tau}}(s),\mathsf{\kappa}^{(0)}(s,t),s)\cdot\mathsf{\xi}^{(0)}(s,t)+\partial_\kappa \mathcal{N}(\overset{\scriptscriptstyle \circ}{\mathsf{\tau}}(s),\mathsf{\kappa}^{(0)}(s,t),s)\cdot\mathsf{\kappa}^{(1)}(s,t)).
\end{align*}
Here, the second term vanishes since $\mathcal{N}(\overset{\scriptscriptstyle \circ}{\mathsf{\tau}}(s),\kappa,s)=\mathsf{0}$ for all $\kappa$ and $s$ (Assumption~\ref{ass:1}).
So, the expression $\mathsf{\xi}=\mu^2\mathsf{L}^{-1}\cdot \mathsf{n}$ with $\mathsf{L}=\partial_\tau\mathcal{N}(\overset{\scriptscriptstyle \circ}{\mathsf{\tau}}(\cdot),\mathsf{\kappa},\cdot)$ is exact up to an error of $\mathcal{O}(\epsilon^2)$. Inserting this expression for $\mathsf{\xi}$ into \eqref{eq:dimless_Rod} yields consequently a $\epsilon$-dependent consistent system up to order $\mathcal{O}(\epsilon^4)$
\begin{equation} \label{eq:eps}
\begin{aligned} \mathsf{D} \cdot \partial_t \mathsf{\breve r} = \mathsf{v},\hspace*{1cm} &\hspace*{4cm}\mathsf{D} \cdot \partial_s \mathsf{\breve r} = \overset{\scriptscriptstyle \circ}{\mathsf{\tau}}+\epsilon^2\mu^2 \mathsf{L}^{-1}\cdot \mathsf{n}\\
\partial_t \mathsf{D} = -\mathsf{\omega}\times \mathsf{D}, \hspace*{1cm} &\hspace*{4cm}  \partial_s \mathsf{D} = -\mathsf{\kappa}\times \mathsf{D}  \\
\epsilon^2\mu^2 \mathsf{L}^{-1}\cdot \partial_t \mathsf{n} &= \partial_s \mathsf{v}+\mathsf{\kappa}\times \mathsf{v} + \overset{\scriptscriptstyle \circ}{\mathsf{\tau}} \times \omega+\epsilon^2\mu^2\mathsf{L}^{-1}\cdot \mathsf{n}\times \omega\\
\partial_t \mathsf{\kappa} &= \partial_s \mathsf{\omega} +\mathsf{\kappa}\times \mathsf{\omega} \\[1.5ex] 
(\rho A)\partial_t\mathsf{v} &= \partial_s \mathsf{n} + \mathsf{\kappa}\times \mathsf{n} +(\rho A) \mathsf{ v} \times \mathsf{\omega}+\mathsf{D}\cdot \mathsf{\breve f}\\  
\epsilon^2(\rho \mathsf{J})\cdot\partial_t \mathsf{\omega} &= \partial_s \mathsf{m} + \mathsf{\kappa}\times \mathsf{m} + \overset{\scriptscriptstyle \circ}{\mathsf{\tau}}\times
 \mathsf{n}+\epsilon^2(\mu^2 \mathsf{L}^{-1}\cdot\mathsf{n}\times \mathsf{n} + (\rho \mathsf{J})\cdot \mathsf{\omega}\times\mathsf{\omega})   
 \end{aligned}
 \end{equation}
with
\begin{align*}
\mathsf{m}(s,t)&=\mu^{-2}\mathcal{M}(\overset{\scriptscriptstyle \circ}{\mathsf{\tau}}(s)+\epsilon^2\mu^2\mathsf{L}^{-1}\cdot\mathsf{n}(s,t),\mathsf{\kappa}(s,t),s), \qquad \qquad 
\mathsf{L}(s,t)=\partial_\tau\mathcal{N}(\overset{\scriptscriptstyle \circ}{\mathsf{\tau}}(s),\mathsf{\kappa}(s,t),s). \nonumber
\end{align*}
Setting $\epsilon = 0$ in \eqref{eq:eps}, we obtain the low-Mach-number--slenderness limit for elastic bodies. The limit equations \eqref{eq:eps0} describe a simplified model of the special Cosserat rod theory, where the angular inertia terms vanish and instead of a material law for the contact force the time-independence of the strain field is imposed ($\partial_t \overset{\scriptscriptstyle \circ}{\mathsf{\tau}}=\mathsf{0}$). This property enforces the limit beam to be inextensible and unshearable and is known as Kirchhoff constraint. Classically, the Kirchhoff constraint relates the tangent and the director triad $\boldsymbol{\tau}=\mathbf{d_3}$, it is often stated as $\mathsf{D} \cdot \partial_s \mathsf{\breve r} = \mathsf{\tau}=\mathsf{e_3}$ in the director basis, see e.g.\ \cite{antman:b:2006, coleman:p:1993, langer:p:1996}. This corresponds to a straight rod curve with perpendicular cross-sections as choice for a stress-free configuration, $\overset{\scriptscriptstyle \circ}{\mathsf{\tau}}=\mathsf{e_3}$.
\begin{equation}\label{eq:eps0}
\begin{aligned} 
\mathsf{D} \cdot \partial_t \mathsf{\breve r} = \mathsf{v},\hspace*{1cm}&\hspace*{4cm} \mathsf{D} \cdot \partial_s \mathsf{\breve r} = \overset{\scriptscriptstyle \circ}{\mathsf{\tau}}\\ 
\partial_t \mathsf{D} = -\mathsf{\omega}\times \mathsf{D}, \hspace*{1cm}&\hspace*{4cm} \partial_s \mathsf{D} = -\mathsf{\kappa}\times \mathsf{D}  \\ 
\mathsf{0} &= \partial_s \mathsf{v}+\mathsf{\kappa}\times \mathsf{v} + \overset{\scriptscriptstyle \circ}{\mathsf{\tau}} \times \omega\\
\partial_t \mathsf{\kappa} &= \partial_s \mathsf{\omega} +\mathsf{\kappa}\times \mathsf{\omega} \\[1.5ex]
(\rho A)\partial_t\mathsf{v} &= \partial_s \mathsf{n} + \mathsf{\kappa}\times \mathsf{n} +(\rho A) \mathsf{ v} \times \mathsf{\omega}+\mathsf{D}\cdot \mathsf{\breve f}\\ 
\mathsf{0} &= \partial_s \mathsf{m} + \mathsf{\kappa}\times \mathsf{m} + \overset{\scriptscriptstyle \circ}{\mathsf{\tau}}\times
 \mathsf{n}  
 \end{aligned}
 \end{equation}
with
\begin{align*}
\mathsf{m}(s,t)&=\mu^{-2}\mathcal{M}(\overset{\scriptscriptstyle \circ}{\mathsf{\tau}}(s),\mathsf{\kappa}(s,t),s).
\end{align*}

In the formulation of \eqref{eq:eps}  the contact force $\mathsf{n}$ and the strain $\mathsf{\tau}$ change their roles, their proportionality operator $\mathsf{L}$ is thereby determined by the material law $\mathcal{N}$. In the limit $\epsilon = 0$ in \eqref{eq:eps0}, $\mathsf{n}$ becomes a variable to the constraint $\partial_t\overset{\scriptscriptstyle \circ}{\mathsf{\tau}}=0$, whereas the material law decouples and is not longer needed for the determination of the solution. The structure of \eqref{eq:eps} obviously changes from a hyperbolic to a degenerate differential-hyperbolic-like system with constraint as $\epsilon \rightarrow 0$, the system becomes stiff.
Note that $\mathsf{L}$ can here be identified with the moduli of the linear elasticity theory, linear operator $(\mathsf{EA})$, cf.\ Example~\ref{ex:1}. But in contrast to the Timoshenko beam whose validity is restricted to small deformations, system \eqref{eq:eps} allows for large deformation as $\mathcal{N}$ can be chosen arbitrarily and $\mathsf{L}$ is just its linearization, $\mathsf{L}=\partial_\tau\mathcal{N}(\overset{\scriptscriptstyle \circ}{\mathsf{\tau}}(\cdot),\cdot)$. Without loss of generality, we might thus restrict our considerations for small $\epsilon$ to an (affine) linear material law of the form $\mathcal{N}(\mathsf{\tau},\mathsf{\kappa},s)=\mathsf{L}(s)\cdot (\mathsf{\tau}-\overset{\scriptscriptstyle \circ}{\mathsf{\tau}}(s))$ which implies an explicit and easily invertible relation between contact force and strain. This relation offers advantages in setting up a numerical scheme that is applicable to both, limit and $\epsilon$-dependent system as $\epsilon\rightarrow 0$.

The presented low-Mach-number--slenderness limit finds its analogue in fluid dynamics where the low-Mach-number limit represents the incompressible limit for compressible fluids  \cite{lions:p:1998,alazard:p:2006}. Here, in elasticity, the asymptotic derivation goes with slow dynamics and slenderness. The limit is valid for arbitrary elastic material laws $\mathcal{M}$ for the contact couple. When the limit system \eqref{eq:eps0} is equipped with the Euler-Bernoulli relation it describes a Kirchhoff beam. In that case with $\overset{\scriptscriptstyle \circ}{\mathsf{\tau}}=\mathsf{e_3}$ the system is also known as Kirchhoff-Love equations \cite{landau:b:1970}. However, in the classical theory the terminology Kirchhoff beam is much wider and stands for an inextensible and unshearable rod (kinetic analogue \cite{antman:b:2006}). The vanishing of the angular inertia terms is generally not presupposed for a Kirchhoff beam (see e.g.\ \cite{dichmann:p:1996, dichmann:p:1996a}), but results here as consequence of the chosen scaling and the associated asymptotics.

%%%%%%%%%%%
\subsection{Asymptotic framework with Euler-Bernoulli material law}

In the further work we focus on the asymptotic framework in a special case, where we restrict on a linear Euler-Bernoulli relation $\mathcal{M}=(\mathsf{EJ})\cdot \mathsf{\kappa}$ and on an homogeneous thread with circular cross-sections and $\overset{\scriptscriptstyle \circ}{\mathsf{\tau}}=\mathsf{e_3}$. Note that these are just technical simplifications to facilitate the numerical studies. They can be easily dropped if it is relevant for certain applications. 
 
In this case the dimensionless quantities become 
\begin{align*}
(\rho A)=1, \quad \quad (\rho \mathsf{J})=\mathsf{P}_2, \quad \quad
\mathsf{L}=(\mathsf{EA})=a^{-1} \mathsf{P}_a, \quad \quad (\mathsf{EJ})= \mathsf{P}_{2/a}.
\end{align*}
where $\mathsf{P}_k=\mathrm{diag}(1,1,k)$, $k\in\mathbb{R}$ and $a=2(1+\nu)$ with Poisson ratio $\nu\in[0,0.5)$ -- in accordance to Example~\ref{ex:1}. The physical quantities stated in Example~\ref{ex:1} correspond to the reference $_\star$-values chosen for the scaling. Hence, we deal with the following system in the (M-T-S) notation that need to be supplemented with appropriate, consistent initial and boundary conditions, it describes the Cosserat rod for $\epsilon>0$ and the Kirchhoff beam for $\epsilon=0$, cf.\ \eqref{eq:eps}, \eqref{eq:eps0} 
\begin{equation}\label{eq:epsEB}
\begin{aligned}
\mathsf{D} \cdot \partial_t \mathsf{\breve r} = \mathsf{v},\hspace*{1cm}&\hspace*{4cm}\mathsf{D} \cdot \partial_s \mathsf{\breve r} = \mathsf{e_3}+\epsilon^2\mu^2a\mathsf{P}_{1/a}\cdot \mathsf{n}\\ 
\partial_t \mathsf{D} = -\mathsf{\omega}\times \mathsf{D}, \hspace*{1cm}&\hspace*{4cm}\partial_s \mathsf{D} = -\mathsf{\kappa}\times \mathsf{D}  \\ 
\epsilon^2 \mu^2a \mathsf{P}_{1/a}\cdot \partial_t \mathsf{n} 
&= \partial_s \mathsf{v}+\mathsf{\kappa}\times \mathsf{v} + \mathsf{e_3} \times \omega + \epsilon^2\mu^2a (\mathsf{P}_{1/a}\cdot \mathsf{n})\times \mathsf{\omega}\\
\partial_t \mathsf{\kappa} &= \partial_s \mathsf{\omega} +\mathsf{\kappa}\times \mathsf{\omega} \\[1.5ex]
\partial_t\mathsf{v} &= \partial_s \mathsf{n} + \mathsf{\kappa}\times \mathsf{n} + \mathsf{ v} \times \mathsf{\omega}+\mathsf{D}\cdot \mathsf{\breve f}\\
\epsilon^2 \mathsf{P}_2\cdot\partial_t \mathsf{\omega} 
&= \mu^{-2}\left(\mathsf{P}_{2/a}\cdot \partial_s \mathsf{\kappa} + \mathsf{\kappa}\times (\mathsf{P}_{2/a}\cdot \mathsf{\kappa})\right) + \mathsf{e_3}\times \mathsf{n}\\
&\quad + \epsilon^2 \left(\mu^2a (\mathsf{P}_{1/a}\cdot\mathsf{n})\times \mathsf{n} + (\mathsf{P}_2\cdot \mathsf{\omega})\times\mathsf{\omega}\right)   
\end{aligned}
\end{equation}

\begin{remark}[Re-formulations of the limit system]\label{rem:limit}
In the limit system, $\epsilon=0$, the Kirchhoff constraint poses a geometric relation between curve and director triad in favor of a material law for the contact force. This motivates the so-called centerline-angle representation of the elastic Kirchhoff beam \cite{langer:p:1996} that renounces the evaluation of the director triad.  Alternatively, the Kirchhoff-Love equations might be known in the invariant form of a wavelike equation for $\mathbf{r}$ with small elliptic regularization due to the bending stiffness (incorporated in $\mu$) \cite{klar:p:2009, marheineke:p:2006}
\begin{align}\label{eq:Fidyst}
\partial_{tt}\mathbf{r}=\partial_s(T\partial_s\mathbf{r})-\mu^{-2}\partial_{ssss}\mathbf{r}+M(\partial_s\mathbf{r} \times \partial_{sss} \mathbf{r})+\mathbf{f}, \qquad  \partial_s M=0, \qquad  \|\partial_s\mathbf{r}\|=1
\end{align}
with torsion couple $M = 2(\mu^2a)^{-1}\,\kappa_3$. The modified traction $T=\mathbf{n}\cdot\mathbf{d_3}-\mu^{-2}\|\partial_{ss}\mathbf{r}\|^2$ acts thereby as Lagrange multiplier to the constraint. Obviously, both re-formulations of the limit system reduce the number of variables, but  strongly change the equations' structure. Appropriate numerical schemes can be found in the respective literature. However, through the low-Mach-number--slenderness limit that clarifies the asymptotic relation between \eqref{eq:R_Cosserat} and \eqref{eq:Fidyst} a uniform numerical treatment for $\epsilon \rightarrow 0$ is made possible. 
\end{remark}

The asymptotic framework provides a special structure of the equations that is exploited for the numerical handling.
Let $\mathsf{\Phi}$ denote the vector-valued function comprising all system variables, $\mathsf{\Phi}(s,t)\in \mathbb{R}^m$, then the $\epsilon$-dependent partial differential algebraic system \eqref{eq:epsEB} can be summarized for all model variants (M), (T), (S) in the form
\begin{align}\label{eq:phi}
\mathsf{A}_\epsilon \cdot \partial_t \mathsf{\Phi}+\mathsf{B}_\epsilon\cdot \partial_s \mathsf{\Phi}+\mathsf{c}_\epsilon(\mathsf{\Phi})=\mathsf{0},
\end{align}
where $\mathsf{A}_\epsilon$ and $\mathsf{B}_\epsilon$ are -- possibly singular -- matrices with constant coefficients and $\mathsf{c}_\epsilon$ is a vector-valued nonlinear function in $\mathsf{\Phi}$. Their $\epsilon$-dependence can be expressed as $\mathsf{A}_\epsilon=\mathsf{A}^{(0)}+\epsilon^2\mathsf{A}^{(1)}$, analogously for $\mathsf{B}_\epsilon$ and $\mathsf{c}_\epsilon$. Expanding $\mathsf{\Phi}$ in the even power series $\mathsf{\Phi}=\sum_{i=0}^n \epsilon^{2i} \mathsf{\Phi}^{(i)}$ and plugging it in \eqref{eq:phi} -- as before --, we find the Kirchhoff beam in $\mathcal{O}(1)$ and its first-order correction in $\mathcal{O}(\epsilon^2)$
 \begin{align}\label{eq:phi0}
 \mathsf{A}^{(0)}\cdot \partial_t \mathsf{\Phi}^{(0)}+\mathsf{B}^{(0)}\cdot \partial_s \mathsf{\Phi}^{(0)}+\mathsf{c}^{(0)}(\mathsf{\Phi}^{(0)})\hspace*{1.3cm}&=\mathsf{0},\\\label{eq:phi1}
 \mathsf{A}^{(0)}\cdot \partial_t \mathsf{\Phi}^{(1)}+\mathsf{B}^{(0)}\cdot \partial_s \mathsf{\Phi}^{(1)}+\partial_{\mathsf \Phi}\mathsf{c}^{(0)}(\mathsf{\Phi}^{(0)})\cdot \mathsf{\Phi}^{(1)}&=\mathsf{f} [\mathsf{\Phi}^{(0)}], \\\nonumber
\text{with} \quad -\mathsf{f}[\mathsf{\Phi}^{0}]&=\mathsf{A}^{(1)}\cdot \partial_t \mathsf{\Phi}^{(0)}+\mathsf{B}^{(1)}\cdot \partial_s \mathsf{\Phi}^{(0)}+\mathsf{c}^{(1)}(\mathsf{\Phi}^{(0)})
 \end{align}
where $\partial_{\mathsf \Phi}\mathsf{c}^{(0)}$ denotes the Jacobi matrix of $\mathsf{c}^{(0)}$. Both asymptotic systems together \eqref{eq:phi0}-\eqref{eq:phi1} are exact up to order $\mathcal{O}(\epsilon^4)$. As usual for asymptotic considerations, the systems have a similar equation structure with the same system matrices  $\mathsf{A}^{(0)}$ and $\mathsf{B}^{(0)}$, they just differ in the right-hand-side and the dependence on the variable. The first-order correction \eqref{eq:phi1} is trivially linear in the variable, whereas the limit system \eqref{eq:phi0} keeps the nonlinearity of the original $\epsilon$-dependent system \eqref{eq:phi}. Having a Newton method for the numerical treatment of the nonlinear term $\mathsf{c}^{(0)}(\mathsf{\Phi}^{(0)})$ in mind, we need its Jacobi matrix already for the computation of the limit system. So, the assembly of the linear system matrix associated with the first-order correction is for free. 

\begin{remark}[Conservation of energy] The Cosserat rod model as well as the Kirchhoff beam equations in \eqref{eq:epsEB} are energy-conserving for external potential forces \cite{simo:p:1988a}, presupposing classical boundary conditions such as, for example, a cantilever beam or a beam with stress-free ends. This holds also true for the system associated with the first-order correction. The corresponding energy with (external) potential energy $V$ is given by $w_\epsilon=w^{(0)}+ \epsilon^2 w^{(1)}+\mathcal{O}(\epsilon^4)$
\begin{align*}
w^{(0)}&=\frac{1}{2}(\mathsf{v}^{(0)})^2+\frac{1}{2\mu^2}\mathsf{\kappa}^{(0)}\cdot\mathsf{P}_{2/a}\cdot \mathsf{\kappa}^{(0)}+V(\mathsf{\breve r}^{(0)})\\
w^{(1)}&=\mathsf{v}^{(0)} \cdot \mathsf{v}^{(1)}
+\frac{1}{\mu^2}\mathsf{\kappa}^{(0)}\cdot\mathsf{P}_{2/a}\cdot \mathsf{\kappa}^{(1)}
+\frac{1}{2}\mathsf{\omega}^{(0)}\cdot\mathsf{P}_{2}\cdot \mathsf{\omega}^{(0)}
+\frac{\mu^2 a}{2} \mathsf{n}^{(0)}\cdot\mathsf{P}_{1/a}\cdot \mathsf{n}^{(0)}
+V(\mathsf{\breve r}^{(1)}).
\end{align*}
\end{remark}

%%%%%%%%%%%%%%%%%%%%%%%%%%%%%%%%%%%%%
\setcounter{equation}{0} \setcounter{figure}{0}
\section{Numerical studies}\label{sec:3}

Being interested in the numerical performance of the asymptotic framework as $\epsilon \rightarrow 0$, we study and discuss the asymptotic convergence and efficiency of the approach in this section. To make use of the structure of the asymptotic systems we apply a Newton method to the discretized equations. The discretization is replaceable. We apply here an asymptotic-preserving scheme in the spirit of the generalized $\alpha$-method \cite{chung:p:1993, sobottka:p:2008} to the underlying systems of partial differential algebraic equations. This scheme leaves the freedom to be switched between energy-conserving and dissipative. This feature can be advantageous when dealing with applications with external non-potential forces (such as fiber-fluid interactions in fiber spinning  \cite{marheineke:p:2006}, non-woven manufacturing \cite{klar:p:2009}, paper making \cite{hamalainen:p:2011} or hair simulations \cite{bertails:p:2006}). For details on the numerical method we refer to the Appendix. As test case we use the two-dimensional Euler-Bernoulli cantilever beam under gravity \cite{fuetterer:p:2012}. This test case offers the advantage that the number of model variables can be reduced from $m=19$ to $9$, while the $\epsilon$-dependent structure of the equations that is of interest is kept. So, clarity is given for the investigation.

\subsection{Test case}\label{sec:test}
\begin{figure}[b]
\vspace*{-0.5cm}
\includegraphics[scale=0.65]{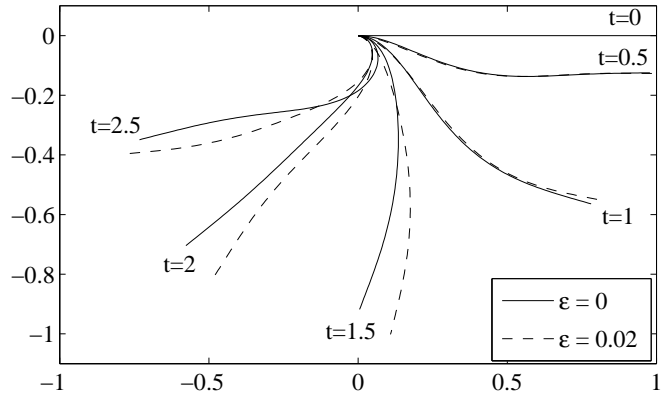}
\vspace*{-1cm}
\caption{\label{fig:faden} Test case: dynamics of the 2d Euler-Bernoulli cantilever beam under gravity $\mathsf{\breve r}$. Simulation of the limit system ($\epsilon=0$) and the $\epsilon$-dependent system ($\epsilon=0.02$).}
\end{figure}

Let $\{\mathbf{e_1},\mathbf{e_2},\mathbf{e_3}\}$ be the outer Cartesian basis. Consider a thread fixed at one end ($s=0$) and stress-free at the other end ($s=1$) that is initially static, stress-free and straight in direction of $\mathbf{e_1}$. Let it be exposed to transversal oscillations in the $\mathbf{e_1}$-$\mathbf{e_2}$-plane due to gravity $\mathbf{f}=-(\rho A)g \mathbf{e_2}$ such that $\mathbf{d_2}=\mathbf{e_3}$ during its motion (cf.\ Figure~\ref{fig:faden}). For this two-dimensional scenario, the model system~\eqref{eq:epsEB} in the (M-T-S) notation simplifies to
\begin{equation*}
\begin{aligned}
\mathsf{D}(\alpha)\cdot \partial_t \mathsf{\breve r} = \mathsf{v},\hspace*{1cm} &\hspace*{4cm}\mathsf{D}(\alpha) \cdot \partial_s \mathsf{\breve r} = \mathsf{e_2}+\epsilon^2\mu^2a\mathsf{P}_{1/a}\cdot \mathsf{n}\\
\partial_t \alpha = \omega, \hspace*{1cm} &\hspace*{4cm}  \partial_s \alpha = \kappa \\
\epsilon^2 \mu^2a \mathsf{P}_{1/a}\cdot\partial_t \mathsf{n} &= \partial_s\mathsf{v}-\kappa\mathsf{v}^\perp-\omega\mathsf{e_1} + \epsilon^2 \mu^2\omega\mathsf{P}_a\cdot\mathsf{n}^\perp \\
\partial_t \kappa &= \partial_s\omega \\[1.5ex]
\partial_t \mathsf{v} &= \partial_s \mathsf{n} - \kappa \mathsf{n}^\perp + \omega \mathsf{v}^\perp + \mathsf{D}(\alpha)\cdot\mathsf{\breve f} \\
\epsilon^2  \partial_t \omega &= \mu^{-2}\partial_s \kappa + n_1+ \epsilon^2 \mu^2 (1-a)n_1n_3 
\end{aligned}
\end{equation*}
with the dimensionless force $\breve{\mathsf{f}}=-\mathrm{Fr}^{-2}\mathsf{e_2}$ and the respective initial and boundary conditions
\begin{align*}
\mathsf{\breve r}(s,0)&=s\mathsf{e_1}, \qquad \quad \alpha(s,0)=\kappa(s,0)=\omega(s,0)=0, \qquad \quad \mathsf{n}(s,0)=\mathsf{v}(s,0)=\mathsf{0}\\
\mathsf{\breve r}(0,t)&=\mathsf{v}(0,t)=\mathsf{0}, \qquad \alpha(0,t)=\omega(0,t)=0,  \qquad \mathsf{n}(1,t)=\mathsf{0}, \qquad  \kappa(1,t)=0.
\end{align*}
Here, we use the abbreviations
$$\mathsf{\breve r}=(\breve r_1,\breve r_2),\quad \mathsf{\breve f}=(\breve f_1,\breve f_2),\quad \mathsf{v}=(v_1,v_3),\quad \mathsf{n}=(n_1,n_3),\quad \omega=\omega_2,\quad \kappa=\kappa_2,\quad m=m_2,$$
$\mathsf{P}_k=\mathrm{diag}(1,k)$, $k=\mathbb{R}$. Moreover, $\mathsf{z}^\perp=(-z_3,z_1)$  denotes the tuple perpendicular to $\mathsf{z}$ for $\mathsf{z}\in\{\mathsf{v},\mathsf{n}\}$. The rotation matrix $\mathsf{D}$ is expressed in terms of the single angle $\alpha$ 
$$\mathsf{D}(\alpha)=\left(\begin{array}{cc} -\sin\alpha & \cos\alpha \\ \cos\alpha &\sin \alpha \end{array}\right), \qquad \alpha=\angle(\mathbf{e_1},\mathbf{d_3}). $$
The Froude number $\mathrm{Fr}$ represents the ratio of inertia and gravity. 

Because the asymptotic results turn out not to depend very sensitively on the specific choice of parameters, we use the setting $(\mathrm{Fr},\mu,a)=(1,10,2.5)$ with the end time $T=2.5$ as example in the following studies. This parameter setting is characteristic in the context of non-woven manufacturing \cite{marheineke:p:2011, klar:p:2009}. Typical properties of polymer fibers are $d=10^{-5}$\,m (diameter), $L=10^{-1}$\,m, $\rho =10^3$\,kg/m$^3$, $E=10^{10}$\,kg/(m\,s$^2$), $\nu=0.25$ and $V=1$\,m/s, implying $\epsilon \sim \mathcal{O}(10^{-4})$. The end time $T$ is chosen respectively with regard to stability issues (see Remark~\ref{rem:oscillations}). The dynamics of the cantilever beam under gravity is illustrated for $\epsilon=0$ and $\epsilon=0.02$ in Figure~\ref{fig:faden}.

\begin{remark}\label{rem:oscillations}
For the planar inextensible beam an existence proof is derived in the absence of gravity and for non-vanishing angular inertia in \cite{caflisch:p:1984}.  Nonlinear planar and non-planar responses of the inextensible beam were investigated numerically in \cite{nayfeh:p:1998} using a combination of the Galerkin procedure and the method of multiple scales, it turned out that the nonlinear geometric terms produce a hardening effect and dominate the non-planar responses for all modes. The inertia terms particularly dominate the high-frequency modes. A global bifurcation analysis for a cantilever beam subjected to a harmonic axial excitation and transverse excitations at the free end was performed in \cite{zhang:p:2005}. Taking into account these results, the simulations of our limit system that has no angular inertia are run in the stable planar regime.
\end{remark}

\subsection{Results}
\begin{figure}[t]
\includegraphics[width=0.45\textwidth]{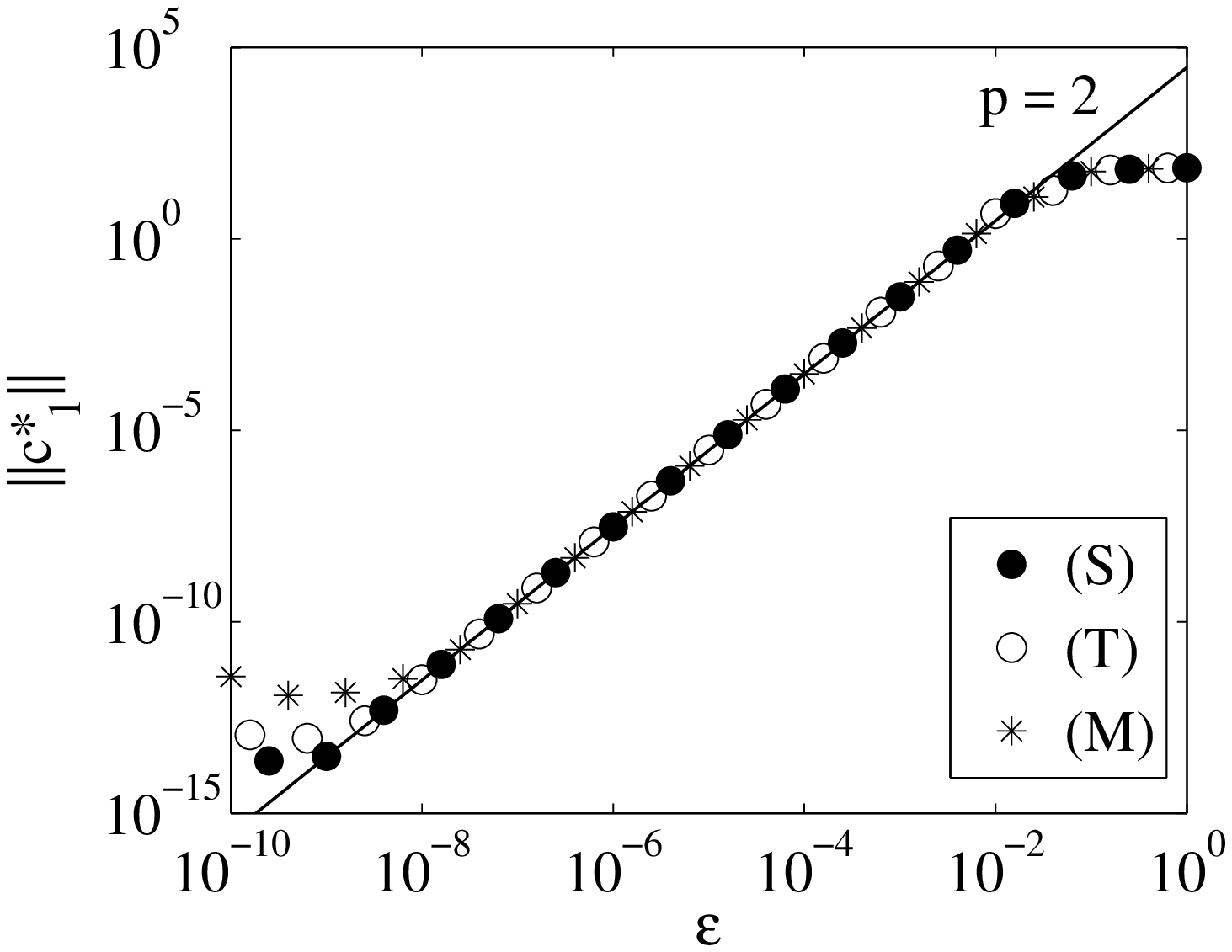}
\includegraphics[width=0.45\textwidth]{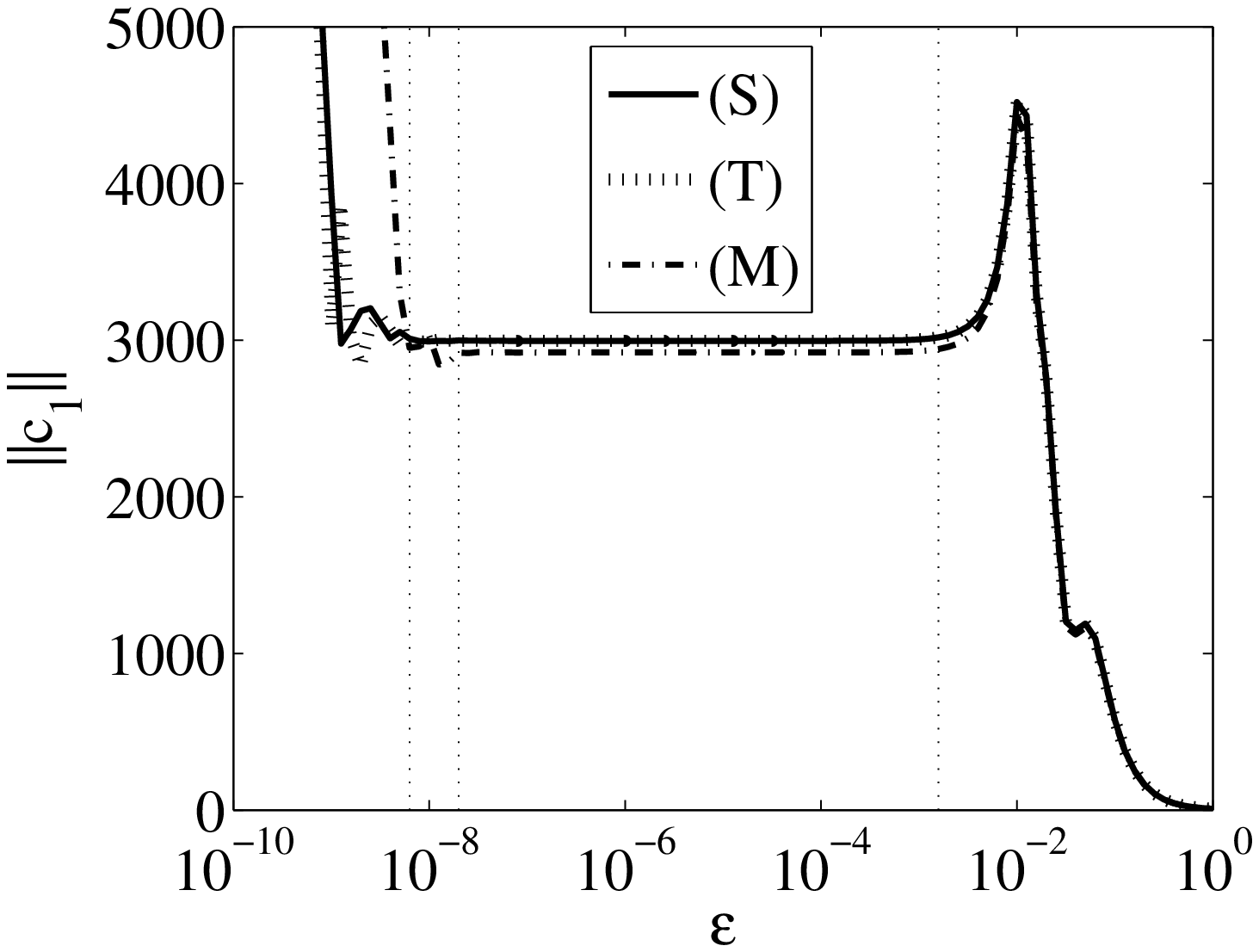}\\
\includegraphics[width=0.45\textwidth]{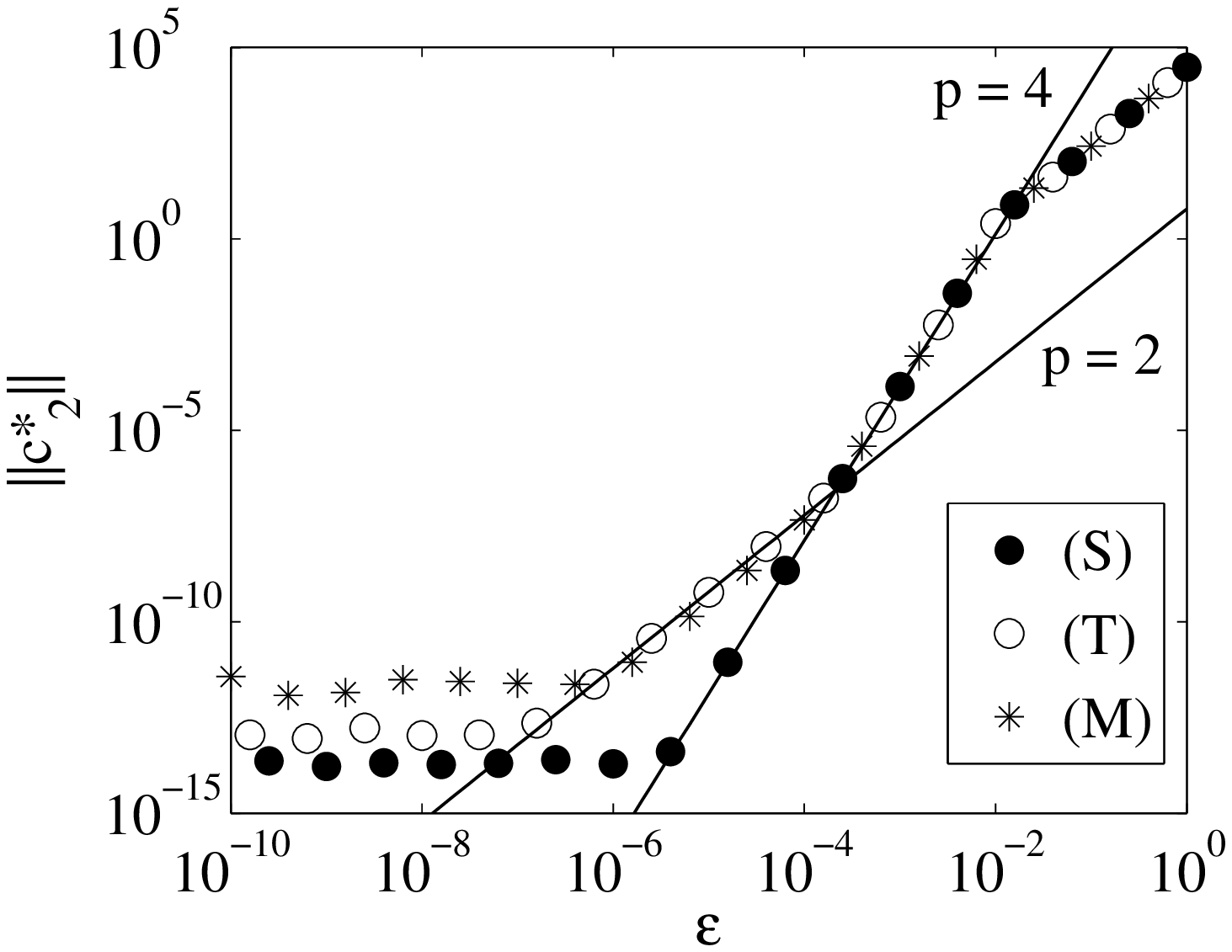}
\includegraphics[width=0.45\textwidth]{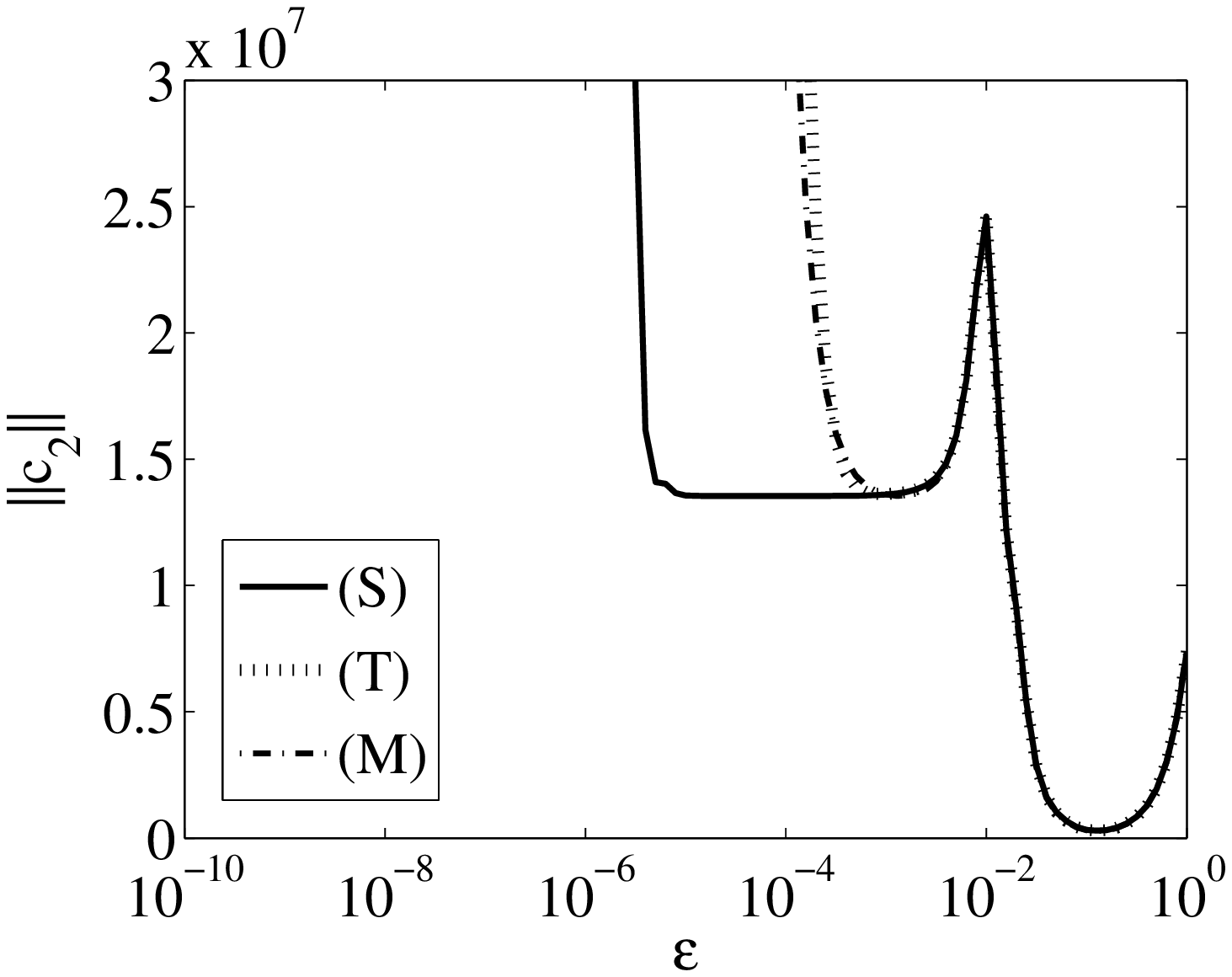}
\caption{\label{fig:eps}Conservation of the asymptotic relations for the different model variants (cf.\ \eqref{eq:c1}-\eqref{eq:c2}, $\mathbf{c_i^\star}=\epsilon^{2i} \mathbf{c_i}$). \emph{Top}: first term $\|\mathbf{c_1^\star}\|$ \emph{(left)} with magnitude $\|\mathbf{c_1}\|$ \emph{(right)} plotted over $\epsilon$. The solid line with $p=2$ indicates quadratic convergence. \emph{Bottom}: second term $\|\mathbf{c_2^\star}\|$ \emph{(left)} with magnitude $\|\mathbf{c_2}\|$ \emph{(right)} plotted over $\epsilon$. The solid lines with $p=4$ and $p=2$ indicate quartic and quadratic convergence, respectively.}
\end{figure}

In the asymptotic framework we have $\mathsf{\Phi}_\epsilon=\mathsf{\Phi}^{(0)}+\epsilon^2\mathsf{\Phi}^{(1)}+\mathcal{O}(\epsilon^4)$, i.e.\ the limit system and its first-order correction are analytically exact up to the order $\mathcal{O}(\epsilon^4)$, cf.\ \eqref{eq:phi}-\eqref{eq:phi1}. The solutions of the limit system $\mathsf{\Phi}^{(0)}$ (Kirchhoff beam), its first-order correction  $\mathsf{\Phi}^{(1)}$ as well as of the original $\epsilon$-dependent system $\mathsf{\Phi}_\epsilon$ (Cosserat rod) are numerically approximated by $\boldsymbol{\varphi}^{(0)}$, $\boldsymbol{\varphi}^{(1)}$ and $\boldsymbol{\varphi}_\epsilon$. We observe that the numerical approximations inherit the asymptotic relation as desired, this means that
\begin{align}\label{eq:c1}
\boldsymbol{\varphi}_\epsilon-\boldsymbol{\varphi}^{(0)}& = \mathbf{c^\star_1} \sim \mathcal{O}(\epsilon^2), 
&& \mathbf{c^\star_1} =\epsilon^2 \, \mathbf{c_1}\\\label{eq:c2}
\boldsymbol{\varphi}_\epsilon-\boldsymbol{\varphi}^{(0)}- \epsilon^2 \boldsymbol{\varphi}^{(1)}&=\mathbf{c^\star_2}\sim \mathcal{O}(\epsilon^4),
&& \mathbf{c^\star_2} = \epsilon^4\,\mathbf{c_2}. 
\end{align}
hold true. Figure~\ref{fig:eps} shows the $\mathcal{L}^2(0,1)$-norm of the vector-valued functions $\mathbf{c^\star_i}$ and $\mathbf{c_i}$, $i=1,2$, in dependence on $\epsilon$ for $\epsilon \in [10^{-10},1]$. The quantities are computed at time $T=2$ for the different model variants (M), (T) and (S) introduced in Notation~\ref{model:rod} by applying the numerical scheme with the parameters $\Delta t= \Delta s=10^{-2}$ and $\lambda=1$. Comparing the model variants, they all yield the quadratic convergence for the first term $\mathbf{c^\star_1}$ with the same ($\epsilon$-independent) magnitude $\|\mathbf{c_1}\|$ in the range $\epsilon \in [10^{-8}, 10^{-2}]$, its relative deviation from $\|\boldsymbol{\varphi}^{(1)}\|$ lies below $10^{-3}$. For the second term $\mathbf{c^\star_2}$ the model variants show the quartic convergence in the range $\epsilon \in [10^{-4},10^{-2}]$ and (S) even in $\epsilon \in [10^{-6},10^{-2}]$. The reduced order of convergence (down to $p=2$ in $\epsilon\in[10^{-6},10^{-4}]$) for  (T) and (M) might be explained by the abandonment of boundary conditions. All boundary conditions are explicitly prescribed and incorporated in the scheme for (S), whereas the boundary conditions associated with $\mathsf{\breve r}$, $\alpha$ for (T) and to $\mathsf{v}$, $\mathsf{n}$ for (M) follow only numerically from the initial conditions and the model equations.  The approximations are consistent but the effect of the term $\epsilon^2 \boldsymbol{\varphi}^{(1)}$ in the numerical differentiation is below the computational accuracy (cancellation error). Also the tails that are observed for very small $\epsilon$ in the convergence plots come from numerical noise. The striking switch in the convergence behavior of both terms $\mathbf{c^\star_i}$, $i=1,2$, that occurs at $\epsilon=10^{-2}$ for all model variants can be explained by the asymptotics itself because the asymptotic framework holds true only for $\epsilon$ sufficiently small. This can be also clearly seen in Figure~\ref{fig:faden} where the beam curve associated with the comparatively large $\epsilon$ lies away from the limit curve. Note that the curves would be indistinguishable for $\epsilon<10^{-2}$. The corresponding peaks observed at $\epsilon=10^{-2}$ in Figure~\ref{fig:eps} (right) leave freedom for speculations: they might be due to the transversal stress component whose effect is then superposed by the other variables.

As the Cosserat rod model contains the small asymptotic parameter explicitly, the system of equations becomes stiff and difficult to solve numerically as $\epsilon\rightarrow 0$. The asymptotic framework provides a special structure of the asymptotic systems \eqref{eq:phi0}-\eqref{eq:phi1} that can be exploited in the numerics when solving the discretized equations with a Newton method. The nonlinear original $\epsilon$-dependent system as well as the nonlinear limit system require 2-3 Newton iterations in average per time step for all model variants, supposing $\Delta t \sim \mathcal{O}(10^{-2})$. The resulting linear systems are of same size and have a similar block-structured band matrix. The Jacobian of the limit system equals the linear system matrix associated with the first-order correction. Therefore, the effort of computing the two asymptotic systems is similar to the effort for the original $\epsilon$-dependent system, as the first-order correction costs only a single linear solving.  Moreover, the solving of the asymptotic systems is much more accurate and robust in determining the influence of small $\epsilon$-values. Consequently, the asymptotic framework makes a uniform numerical handling of the transition regime for small $\epsilon$ easily possible. The choice of the model variant has only a very marginal effect on the numerical performance. It influences neither efficiency nor convergence properties, for details see the Appendix. However, (T) seems to be preferable for energy-conservation and (S) for robustness and asymptotic-conservation.

%%%%%%%%%%%%%%%%%%%%%%%%%%%%%%%%%%%%%%%%%%%%%%%%%%%%%%

\section{Conclusion}
The low-Mach-number--slenderness limit ($\epsilon \rightarrow 0$ and $\mathrm{Ma} \rightarrow 0$ with $\mathrm{Ma}/\epsilon =\mu =\mathrm{const}>0$) describes the asymptotic relation between the dynamic elastic Cosserat rod and the inextensible, unshearable Kirchhoff beam without angular inertia. Its derivation requires the unique determination of the family of stress-free configurations by a strain field and the monotonicity of the constitutive law for the contact force. The limit is valid for general material laws for the contact couple. As the Cosserat rod model contains the small asymptotic parameter explicitly, the system of equations becomes stiff and difficult to solve numerically as $\epsilon\rightarrow 0$. The derived asymptotic framework allows a uniform numerical treatment of the transition regime between the two models. We point out that solving the asymptotic systems (limit system and its first-order correction) that are exact up to order $\mathcal{O}(\epsilon^4)$ is of same computational costs as solving the original $\epsilon$-dependent system, but much more robust since the equations are independent of $\epsilon$. Hence, we propose to use the asymptotic systems for the numerical investigation of applications where the beam velocity compared to the typical speed of sound is small and the slenderness ratio is small, as it is the case, for example, in fiber dynamic simulations of non-woven production processes \cite{klar:p:2009, marheineke:p:2011} or hair modeling in computer graphics \cite{audoly:b:2010, bergou:p:2008}.

%%%%%%%%%%%%%%%%%%%%%%%%%%%%%%%%%%%%%%%%%%%%%%%%%%%%%%

\appendix
\renewcommand{\theequation}{A.\arabic{equation}}
\renewcommand{\thetable}{A.\arabic{table}}
\renewcommand{\thefigure}{A.\arabic{figure}}
\setcounter{equation}{0} \setcounter{figure}{0}

\section{Numerical Method}
The rod models consist of systems of partial differential algebraic equations. The model variants (M), (T) and (S) introduced in Notation~\ref{model:rod} imply temporal or spatial differential-algebraic systems of different index (index 0 to 3) after semi-discretization in space or time, respectively. We propose a numerical scheme that is applicable to all formulations in the asymptotic framework. For the numerical treatment of the respective model equations we combine a Gauss-Legendre collocation method (finite differences) on a equidistant space grid with a flexible time integration. It can be also viewed as a conservative finite-volume method. Using a temporal tuning parameter $\lambda\in [0.5,1]$ the time integration can be switched continuously from an energy-conserving Gauss midpoint rule ($\lambda=0.5$) to a dissipative implicit Euler method ($\lambda=1$) in the spirit of a generalized $\alpha$-method \cite{chung:p:1993, sobottka:p:2008}.  The discretized Cosserat rod and Kirchhoff beam models result in nonlinear systems of equations that are solved with a Newton method with Armijo step size control. 

\begin{remark}
To exploit the structure of the asymptotic systems only the application of the Newton method is essential, whereas the underlying discretization is replaceable.
Hence, we refer to the well-established results by Simo \& Vu-Quoc, Simo et al.\ \cite{simo:p:1986, simo:p:1988, simo:p:1995} or Romero \& Armero \cite{romero:p:2002} for exact energy-momentum conserving algorithms. Preserving conservation properties and objectivity is also the topic in, among others, \cite{betsch:p:2003, leyendecker:p:2006}, where the beam equations are formulated as a constraint Hamilitonian system. The ingredients of the scheme are finite elements and a G-equivariant discrete derivative, the constraints are treated by the Lagrange multiplier method, the penalty method or the augmented Lagrange method. In \cite{lang:p:2011} the Cosserat rod model is transformed into a Lagrangian differential-algebraic system of index 3 that is reduced to index 0 by introducing Baumgarte penalty accelerations for stabilization and solved by using finite differences on a staggered grid. The stiffness of the inextensible, unshearable beam with angular inertia is handled by the schemes in e.g.\ \cite{hou:p:1998,dichmann:p:1996a}. Further approaches based on spatial semi-discretizations and Lagrangian mechanics can be found in literature; for algorithms used in the computer graphics see e.g.\ \cite{bergou:p:2008, bertails:p:2006, spillmann:p:2007}.  
\end{remark}

\subsection{Numerical scheme}
Consider the spatial and temporal grids $s_i=i \Delta s$, $i=0,...,N$ and $t^j=j\Delta t$, $j=0,1,...$ with fixed cell and step size $\Delta s$ and $\Delta t$. Let the vector of the unknown system variables at the grid point $(s_i,t^j)$ be denoted by $\boldsymbol{\varphi}_i^j\in \mathbb{R}^m$ and all spatial information at the time level $t^j$ be summarized in $\boldsymbol{\varphi}^j\in \mathbb{R}^{m(N+1)}$. We treat $(s_{i+1/2},t^{j+\lambda})$ with $\lambda \in [0.5,1]$ as collocation points. The idea of the scheme is to fulfill the partial differential algebraic system $\mathsf{F}(\partial_t \mathsf{\Phi}, \partial_s \mathsf{\Phi}, \mathsf{\Phi})=\mathsf{0}$ at all collocation points and approximate the appearing variables and derivatives at the collocation points in terms of $\boldsymbol{\varphi}_i^j$, $i=0,...,N$ and $j=0,1,...$. 
We have
\begin{align*}
\mathsf{F}((\partial_t \mathsf{\Phi}, \partial_s \mathsf{\Phi}, \mathsf{\Phi})_{i+1/2}^{j+\lambda})&=\mathsf{0}, \quad i=0,...,N-1, \quad j=0,1,...
\end{align*}
where $\mathsf{\Phi}_{i+1/2}^{j+\lambda}$, $\partial_t \mathsf{\Phi}_{i+1/2}^{j+\lambda}$ and $\partial_s\mathsf{\Phi}_{i+1/2}^{j+\lambda}$ are discretized with the following four-point stencils 
\begin{align*}
\mathsf{\Phi}_{i+1/2}^{j+\lambda}&=\frac{\lambda}{2}(\boldsymbol{\varphi}_{i+1}^{j+1}+\boldsymbol{\varphi}_i^{j+1})+\frac{1-\lambda}{2}(\boldsymbol{\varphi}_{i+1}^{j}+\boldsymbol{\varphi}_i^{j})\\
\partial_t \mathsf{\Phi}_{i+1/2}^{j+\lambda}&=\frac{1}{2\Delta t}(\boldsymbol{\varphi}_{i+1}^{j+1}+\boldsymbol{\varphi}_i^{j+1}-\boldsymbol{\varphi}_{i+1}^{j}-\boldsymbol{\varphi}_i^{j})\\
\partial_s\mathsf{\Phi}_{i+1/2}^{j+\lambda}&=\frac{\lambda}{\Delta s}(\boldsymbol{\varphi}_{i+1}^{j+1}-\boldsymbol{\varphi}_i^{j+1})+\frac{1-\lambda}{\Delta s}(\boldsymbol{\varphi}_{i+1}^{j}-\boldsymbol{\varphi}_i^{j}).
\end{align*} 
The stencils are based on the convex approximation $\mathsf{\Phi}_{i}^{j+\lambda}=\lambda \boldsymbol{\varphi}_{i}^{j+1}+(1-\lambda)\boldsymbol{\varphi}_{i}^j$ (analogously for $\mathsf{\Phi}_{i+1/2}^{j}$) and finite differences for the derivative $\partial_t \mathsf{\Phi}_{i}^{j+\lambda}=(\boldsymbol{\varphi}_{i}^{j+1}-\boldsymbol{\varphi}_{i}^j)/\Delta t$ (analogously for $\partial_s \mathsf{\Phi}_{i+1/2}^{j}$). This discretization can be in particular viewed as a conservative finite-volume scheme. 
We obtain the unknowns $\boldsymbol{\varphi}^{j+1}$ on the new time level for given $\boldsymbol{\varphi}^j$ by solving the system
\begin{align}\label{eq:F}
&\mathcal{F}(\boldsymbol{\varphi}^{j+1})=(\mathcal{F}_0(\boldsymbol{\varphi}_0^{j+1}\boldsymbol{\varphi}_1^{j+1}), ...,  \mathcal{F}_{N-1}(\boldsymbol{\varphi}_{N-1}^{j+1}\boldsymbol{\varphi}_N^{j+1}), \mathcal{G}(\boldsymbol{\varphi}_{0}^{j+1}\boldsymbol{\varphi}_N^{j+1}))= \mathsf{0}\\\nonumber
&\text{with}\quad \quad \mathcal{F}_i(\boldsymbol{\varphi}_i^{j+1}\boldsymbol{\varphi}_{i+1}^{j+1})=\mathsf{F}((\partial_t \mathsf{\Phi}, \partial_s \mathsf{\Phi}, \mathsf{\Phi})_{i+1/2}^{j+\lambda}(\boldsymbol{\varphi}^{j+1}_i, \boldsymbol{\varphi}^{j+1}_{i+1})), \quad i=0,...,N-1.
\end{align}
Note that in the definition of $\mathcal{F}_i$ we treat $(\partial_t \mathsf{\Phi}, \partial_s \mathsf{\Phi}, \mathsf{\Phi})_{i+1/2}^{j+\lambda}$ as function in $(\boldsymbol{\varphi}^{j+1}_i, \boldsymbol{\varphi}^{j+1}_{i+1})$ according to the discretization stencils stated above. Moreover, the function $\mathcal{G}$ prescribes the boundary conditions and  $\boldsymbol{\varphi}^0$ contains the initial conditions.

\begin{figure}[t]
\includegraphics[width=0.525\textwidth]{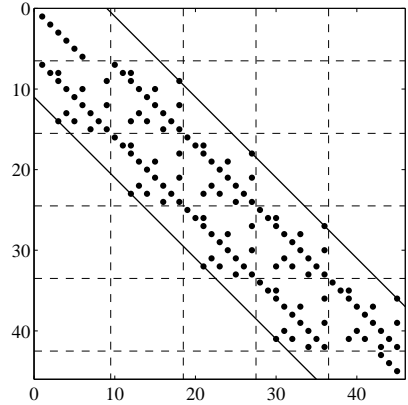}
\vspace*{-0.5cm}
\caption{\label{fig:jacobian}Pattern of the Jacobian for the limit system with $m=9$, $N=4$ and (S) that corresponds to the test case of the 2d Euler-Bernoulli cantilever beam in Section~\ref{sec:test}.}
\end{figure}

In case of the Cosserat rod and Kirchhoff beam models, the function $\mathcal{F}$ is nonlinear and \eqref{eq:F} is solved with a Newton method. The Jacobi matrix $\partial_{\boldsymbol{\varphi}}\mathcal{F}\in \mathbb{R}^{m(N+1)\times m(N+1)}$ is a band matrix with block structure, presupposing an appropriate sorting of the boundary conditions. The blocks are given by  
\begin{align*}
\partial_{\boldsymbol{\varphi}_i} \mathcal{F}_i &=\frac{1}{2\Delta t} \mathsf{A}-\frac{\lambda}{\Delta s}\mathsf{B}+\frac{\lambda}{2} \partial_\mathsf{\Phi}\mathsf{c}(\mathsf{\Phi}_{i+1/2}^{j+\lambda}), \quad \quad \quad 
\partial_{\boldsymbol{\varphi}_{i+1}} \mathcal{F}_i =\frac{1}{2\Delta t} \mathsf{A}+\frac{\lambda}{\Delta s}\mathsf{B}+\frac{\lambda}{2} \partial_\mathsf{\Phi}\mathsf{c}(\mathsf{\Phi}_{i+1/2}^{j+\lambda})
\end{align*}
where $\mathsf{A}$, $\mathsf{B}$, $\partial_\mathsf{\Phi}\mathsf{c}\in \mathbb{R}^{m\times m}$ stand for the system matrices of the $\epsilon$-dependent model and the limit model ($\epsilon=0$), respectively (cf.\ \eqref{eq:phi}-\eqref{eq:phi1}). The blocks associated with the limit model are less occupied since $\mathsf{A}_\epsilon=\mathsf{A}^{(0)}+\epsilon^2 \mathsf{A}^{(1)}$ (analogously for $\mathsf{B}_\epsilon$, $\mathsf{c}_\epsilon$).  Moreover, the Jacobian of the discretized limit model equals the linear system matrix associated with the discretized model of the first-order correction. See Figure~\ref{fig:jacobian} for an example pattern. 

\begin{remark}\label{rem:SD}
The presented approach can be also interpreted as spatial semi-discretization with an initial value problem as well as temporal semi-discretization with a boundary value problem. In any interpretation a differential-algebraic system (DAE) is obtained. Its index varies from 0 to 3, depending on the underlying system ($\epsilon>0$ (Cosserat rod), $\epsilon=0$ (Kirchhoff beam), first-order correction $\mathcal{O}(\epsilon^2)$), the model variant (M), (T), (S) and the considered semi-discretization, see Table~\ref{tab:1}. In this context the used discretization schemes can be understood as implicit Runge-Kutta methods of stage $s=1$: Gauss method for the DAE in space and in time with $\lambda=0.5$ and Radau IIa method (Euler method) for the DAE in time with $\lambda=1$. Although theoretical convergence results are only available up to index 1 for the Gauss method (convergence order $p=2$) and up to index 2 for the Radau IIa method (convergence order $p=1$) \cite{hairer:b:1989}, we point out that in practice the numerical scheme can be successfully applied to all cases.
\end{remark}

\begin{table}[t]\label{tab:1}
\caption{DAE index for the semi-discretized systems (cf.~Remark~\ref{rem:SD})}
\begin{tabular}{|l || c c c | c c c|}
\hline
& \multicolumn{3}{c|}{DAE in time} &  \multicolumn{3}{c|}{DAE in space}\\
& (M) & (T) & (S) & (M) & (T) & (S)\\
\hline \hline
$\epsilon>0$ & 1 & 0 & 1 & 1 & 1 & 0\\
$\epsilon=0$ & 3 & 2 & 2 & 1 & 1 & 0\\
$\mathcal{O}(\epsilon^2)$ & 3 & 2 & 2 & 1 &1 & 0\\
\hline
\end{tabular}
\end{table} 

\subsection{Convergence properties}
\begin{figure}[t]
\includegraphics[width=0.45\textwidth]{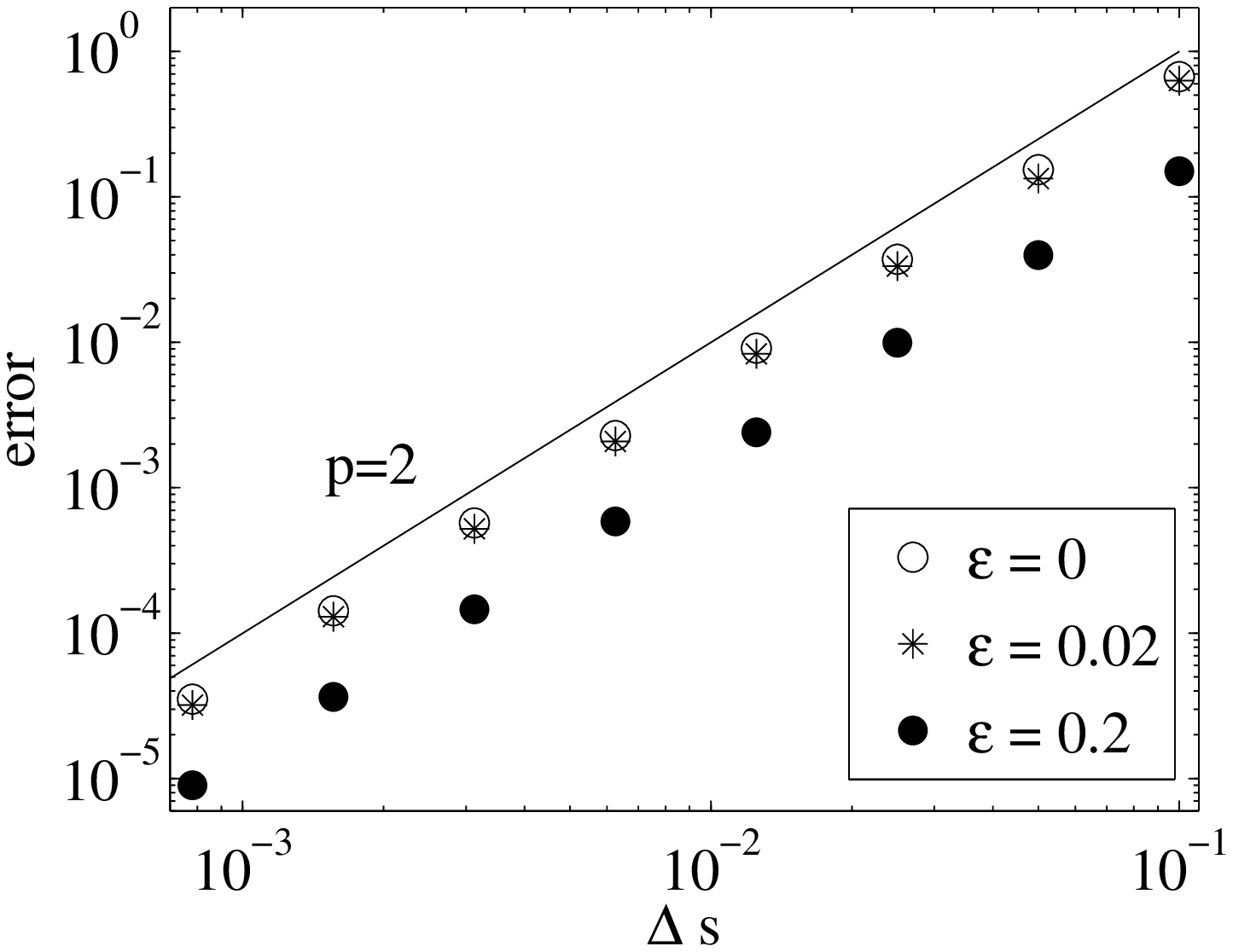}\\
\vspace*{-0.75cm}
\includegraphics[width=0.45\textwidth]{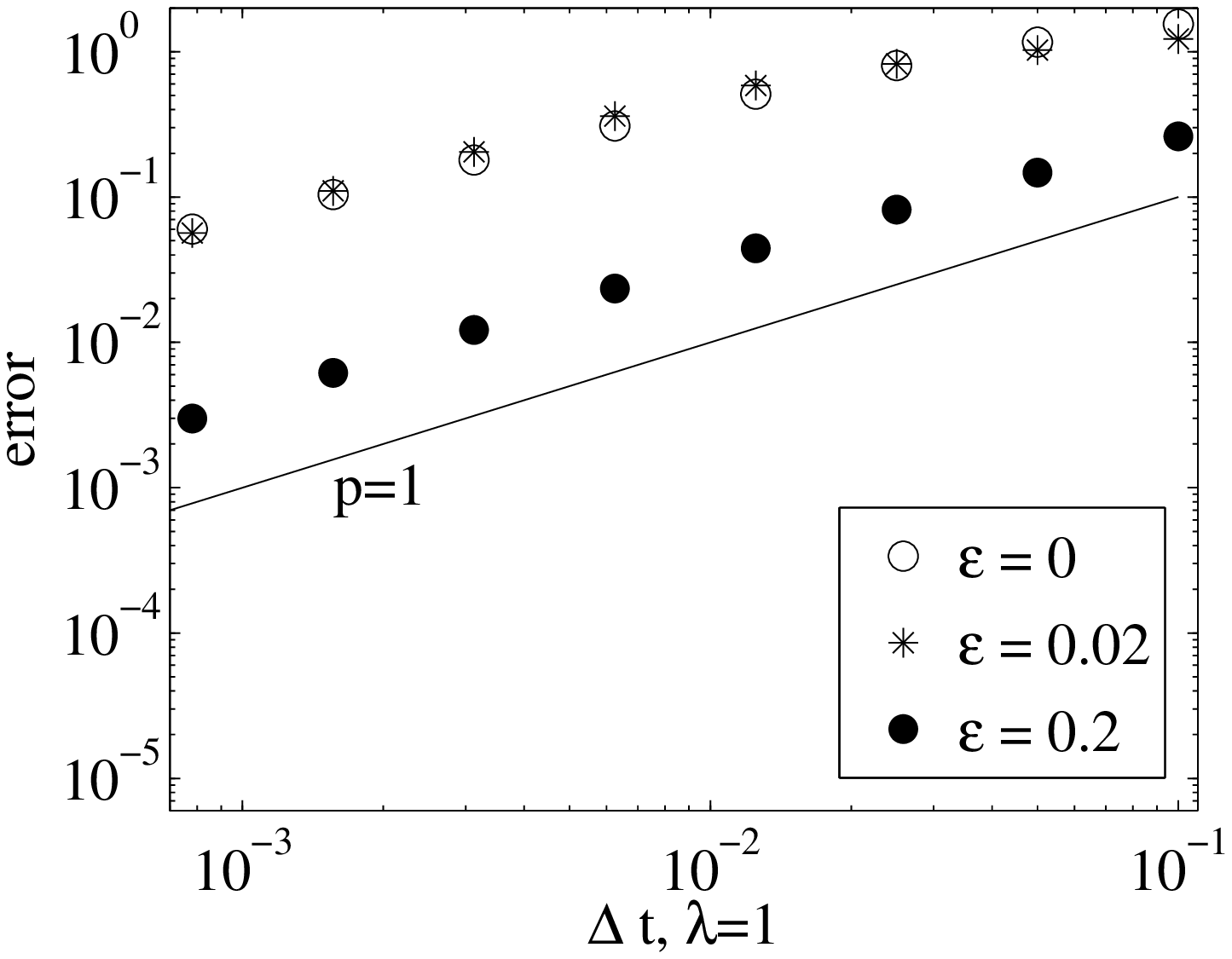}
\includegraphics[width=0.45\textwidth]{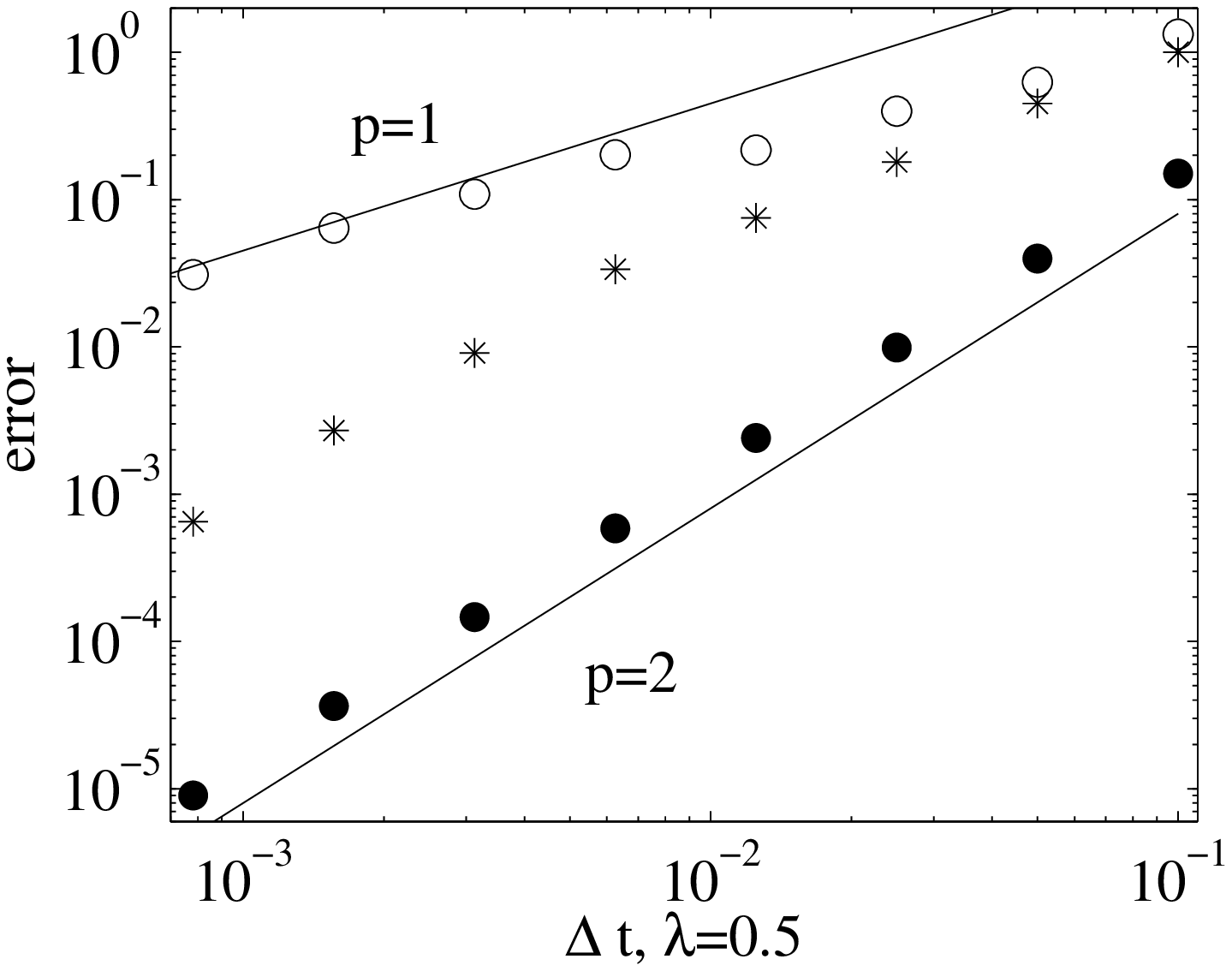}
\vspace*{-1cm}
\caption{Convergence results for the numerical scheme, relative $\mathcal{L}^2(0,1)$-error at time $T$ for $\epsilon$-dependent and limit systems. \emph{Top}: spatial convergence ($\Delta s \rightarrow 0$, fixed $\Delta t$). \emph{Bottom}: temporal convergence ($\Delta t \rightarrow 0$, fixed $\Delta s$) for $\lambda=1$ \emph{(left)} and $\lambda=0.5$ \emph{(right)}. \label{fig:conv}}
\end{figure}

The scheme has the freedom that it can be switched between energy-conserving and dissipative by changing the time integration between the Gauss midpoint rule ($\lambda=0.5$) and the implicit Euler method ($\lambda=1$). This change effects the convergence properties, as we will show exemplarily for the test case in Section~\ref{sec:test}.  We study at first separately the spatial and temporal convergence orders in the context of semi-discretized DAEs (see Remark~\ref{rem:SD}). As expected and in agreement with the analytical results, the midpoint discretization for the spatial DAE yields a numerical convergence of order $p_s=2$ for all variables in both, $\epsilon$-dependent and limit systems. In particular, the results are independent of the chosen model variant (M), (T), (S).  Figure~\ref{fig:conv} (top) shows the relative $\mathcal{L}^2(0,1)$-error between the reference solution associated with $\Delta s_{\mathrm{ref}}=10^{-4}$ and the approximations for $\Delta s \in \{10^{-1}, 0.5^1\cdot10^{-1}, ..., 0.5^7\cdot10^{-1}\}$ at $T=2$ computed with $\Delta t=10^{-2}$, (S) and $\lambda=1$.  In contrast, the temporal convergence studies are more sophisticated. For the $\epsilon$-dependent system the temporal DAE is of index $0$ (T) or $1$ (M), (S). In agreement with the theory we obtain a numerical convergence of order $p_t=2$ for $\lambda=0.5$ and $p_t=1$ for $\lambda=1$. For the limit system the temporal DAE is of index $2$ (T), (S) or even $3$ (M). In spite of the lack of convergence theory we find here $p_t\approx 1$, $p_t<1$. The algebraic variables behave in general worse than the differential ones. Interesting to note is that the results are similar for all model variants. Figure~\ref{fig:conv} (bottom) shows the relative $\mathcal{L}^2(0,1)$-error between the reference solution associated with $\Delta t_{\mathrm{ref}}=10^{-4}$ and the approximations for $\Delta t \in \{10^{-1}, 0.5^1\cdot10^{-1}, ..., 0.5^7\cdot10^{-1}\}$ at $T=2$ computed with $\Delta s=2\cdot 10^{-3}$, (S) and both $\lambda=0.5$ and $\lambda= 1$. On first glance the error values (in particular for the limit system) maybe look large, but this is only due to the used time $T$. Moreover, we point out that the temporal convergence requires a sufficiently small $\Delta s$. Combining temporal and spatial convergence and considering $\Delta t= \Delta s\rightarrow 0$, we always obtain the convergence order $p=1$ for the whole scheme. The only exception is the $\epsilon$-dependent case with moderate $\epsilon$ values solved with $\lambda=0.5$, here we have $p=2$. This coincides with the separate investigations as we get $p=\min\{p_t, p_s\}$. 

%%%%%%%%%%%%%%%%%%%

\subsection*{Acknowledgements}
The authors acknowledge the support by the German BMBF (Project OPAL 05M13) and DFG (Project MA 4526/2-1, WE 2003/4-1).

%%%%
%\bibliographystyle{siam}
%\bibliography{ref_2013}

%%%%

\end{document}